\documentclass{amsart}
\usepackage{amsfonts,amssymb,amsmath,a4wide}
\newcommand{\qedbox}{ \fbox{}}

\newtheorem{teo}{Theorem}[section]
\newtheorem{lema}{Lemma}[section]
\newtheorem{pro}{Proposition}[section]
\newtheorem{defi}{Definition}[section]\newtheorem{rema}{Remark}[section]
\newtheorem{coro}{Corollary}[section]
\numberwithin{equation}{section}

\def \b {\bullet}

\def \V{\mathcal{V}}

\def \bnabla{\overline{\nabla}}

\def \z{\zeta}

\def\J{\mathcal{J}}

\def \l4{[[\Lambda^{0,4}]]}
\def \ll3{[[\Lambda^{0,3}]]}
\def \ho{\otimes_1}
\def \ko{\otimes_2}
\def \bJ{\mathbb{J}}

\def \pf{\partial}
\def \bpf{\overline{\partial}}
\def \g{\mathfrak{g}}
\def \K{\mathcal{K}}
\def \RR {\mathcal{R}}
\def\sideremark#1{\ifvmode\leqslantavevmode\fi\vadjust{\vbox to0pt{\vss% the remark
 \hbox to 0pt{\hskip\hsize\hskip1em%                          will appear only
 \vbox{\hsize2.5cm\tiny\raggedright\pretolerance10000%          on the side
 \noindent #1\hfill}\hss}\vbox to8pt{\vfil}\vss}}}%
\begin{document}\larger[2]
\title[]{Connexions with totally skew-symmetric torsion and nearly-K\"ahler geometry}
\begin{abstract} We study almost Hermitian structures admitting a Hermitian connexion with totally skew-symmetric torsion or equivalently, those
almost Hermitian structures with totally skew-symmetric Nijenhuis tensor. We investigate up to what extent the Nijenhuis tensor fails 
to be parallel with respect to the characteristic connexion. This is naturally described by means of 
an extension of the notion of Killing form to almost Hermitian geometry. In this context, we also make an essentially self-contained survey of nearly-K\"ahler geometry, but from the perspective of non-integrable holonomy systems.
\end{abstract}
\subjclass[2000]{53C12, 53C24, 53C55}
\keywords{nearly K\"ahler manifold, connexion with torsion}
\author[P.- A. Nagy]{Paul-Andi Nagy}
\address[P.-A. Nagy]{Department of Mathematics, University of Auckland, Private Bag 92019, Auckland, New Zealand}
\email{nagy@math.auckland.ac.nz}
\date{\today}
\maketitle
\tableofcontents
\section{Introduction}
For a given Riemannian manifold $(M^n,g)$, the holonomy group $Hol(M,g)$of the metric $g$ contains fundamental geometric information. If the manifold is locally irreducible, which can always be assumed by the well-known 
deRham splitting theorem, assuming that $Hol(M,g)$ is smaller that $SO(n)$ has strong geometric implications at the level of the geometric structure supported by $M$ (see \cite{redbook} for an account of the classification of 
irreducible Riemannian holonomies). Nowadays, this has been extended to cover the classification of possible holonomy groups of torsionless affine connections  (see \cite{bryant1, bryant2, MerkSwach} for an overview). 

From a somewhat different perspective, in the past years there has been increased interest in the theory of metric connexions with totally skew-symmetric torsion on Riemannian manifolds. 
The holonomy group of such a connexion is once again an insightful object which can be used to describe the geometry in many situations.
One general reason supporting this fact is that the deRham splitting theorem fails to hold in this generalised setting and the counterexamples to this extent are more intricate than just Riemannian products.  Additional motivation for the study 
of connections in the above mentioned class comes from string theory, where this requirement is part of various models (see \cite{AFNP, Stro} for instance).

Actually, some of these directions go back to the work of A.Gray, in the seventies, in connection to the so called concept of weak holonomy \cite{g0, Alw}, which was aimed to generalise 
that of Riemannian holonomy. In this respect, the study of nearly-K\"ahler manifolds initiated by Gray in \cite{g1,g2,g3} has been instrumental in the theory of geometric structures of non-integrable type. 
As it is well known, for almost Hermitian structures a classification has been first given by Gray and Hervella in \cite{GrayH}, which proved to be quite 
suitable for generalisation. 

Recently, general $G$-structures of non-integrable type have been studied in \cite{swann, FriIva, cs, Frip}, by taking into account the algebraic properties of the various torsion modules. 
By using results from holonomy theory they have indicated how to read off various intersections of $G$-modules interesting geometric properties  
when in presence, say, of a connexion with totally skew-symmetric torsion. Aspects of such an approach are contained in the Gray-Hervella classification and its analogues for other classical groups, we shall comment on later. Note that a central place in \cite{cs} is occupied by connections with parallel skew-symmetric torsion, studied since then by many authors.

By combining these two points of view, let us return to the Gray-Hervella classification and recall that one of its subclasses, that of nearly-K\"ahler structures, has a priori 
parallel torsion for a connection with totally skew-symmetric torsion \cite{kiri}. This suggests perhaps that geometric properties of the intrinsic torsion might by obtained in more general context.

In this paper we will at first investigate some of the properties of almost Hermitian structures in the Gray-Hervella class $W_1+W_3+W_4$, also referred to -as we shall in what follows-
as the class $\mathcal{G}_1$. It is best described \cite{FriIva} as 
the class of almost Hermitian structures admitting a metric and Hermitian connection with totally  skew-symmetric torsion.  Yet, this property is equivalent \cite{FriIva} with requiring the Nijenhuis tensor of the almost 
complex structure to be a $3$-form when evaluated using the Riemannian metric of the structure. There are several subclasses of interest, for example the class $W_3+W_4$ consisting of Hermitian structures, the class 
$W_4$ of locally conformally K\"ahler metrics and the class $W_1$ of nearly-K\"ahler structures, which will occupy us in the second half of the paper. One general question which arises is how  far the class 
$W_1+W_3+W_4$ is, for instance, from the subclass $W_3+W_4$ of Hermitian structures. 

The paper is organised as follows. In Section 2 we present a number of algebraic facts related to representations of the unitary group and also some background material 
on algebraic curvature tensors. We also briefly review some basic facts from almost-Hermitian geometry, including a short presentation of the Hermitian connections of relevance for 
what follows. At the end of the Section 2 we introduce the notion of Hermitian Killing form with respect to a Hermitian connection with torsion. Although this 
is a straightforward variation on the Riemannian notion it will serve an explanatory r\^ole in the next section. 

Section 3 is devoted to the study of the properties of the torsion tensor of an almost Hermitian structure in the class $\mathcal{G}_1$. More precisely we prove:
\begin{teo} \label{intromain}Let $(M^{2m},g,J)$ be an almost-Hermitian structure of type $\mathcal{G}_1$ and let $D$ be the unique Hermitian connection whose torsion tensor 
$T^D$ belongs to $\Lambda^3(M)$. Then:
$$D_X \psi^{+}=-\frac{1}{2}(X \lrcorner dT^D)_{\lambda^3}
$$
for all $X$ in $TM$, where $\psi^{+}$ is proportional to the Nijenhuis form of $J$ and the subscript indicated orthogonal projection on $\lambda^3(M)$, the bundle of real valued forms of type $(0,3)+(3,0)$. 
\end{teo}
For unexplained notation and the various numerical conventions we refer the reader to Section 2.
In the language of that section, this can be rephrased to say that $\psi^{+}$ is a Hermitian Killing form of type $(3,0)+(0,3)$. Theorem \ref{intromain} is extending to the wider class of $\mathcal{G}_1$-structures the well known \cite{kiri}
parallelism of the Nijenhuis tensor of a nearly-K\"ahler structure. We also give necessary and sufficient conditions for the parallelism of the Nijenhuis tensor, in a $\mathcal{G}_1$-context and specialise Theorem \ref{intromain} to the subclass of 
$W_1+W_4$ manifolds which has been recently subject of attention \cite{But2, CsIv2, MorOr} in dimension $6$. Section 3 ends with the explicit computation of the most relevant components of the curvature 
tensor of the connection $D$ above, in a ready to use form. 

Section 4 contains a survey of available classification results in nearly-K\"ahler geometry. This is mainly based on \cite{Nagy1, Nagy2}, but here we have adopted a slightly different approach, which emphasises the  
r\^ole of the non-Riemannian holonomy system which actually governs the way various holonomy reductions are performed geometrically. We have outlined, step by step, the procedure reducing 
a holonomy system of nearly-K\"ahler type to an embedded irreducible Riemannian holonomy system, which is a the key argument in the classification result which is stated below.
\begin{teo}\cite{Nagy2} \label{classnk}Let $(M^{2m},g,J)$ be a complete SNK-manifold. Then $M$ is, up to finite cover, a Riemannian product whose factors belong to the following classes:
\begin{itemize}
\item[(i)] homogeneous SNK-manifolds;
\item[(ii)] twistor spaces over positive quaternionic K\"ahler manifolds;
\item[(iii)] $6$-dimensional $SNK$-manifolds.
\end{itemize}
\end{teo}
We have attempted to make the text as self-contained as possible, for the convenience of the reader.

\section{Almost Hermitian geometry}
An almost Hermitian structure $(g,J)$ on a manifold $M$ consists in a Riemannian metric $g$ on $M$ and a compatible almost complex structure 
$J$. That is, $J$ is an endomorphism of the tangent bundle to $M$ such that $J^2=-1_{TM}$ and 
$$g(JX,JY)=g(X,Y)$$ 
for all $X,Y$ in $TM$. It follows that the dimension of $M$ is even, to be denoted by $2m$ in the subsequent. Any almost Hermitian structure 
comes equipped with its so-called K\"ahler form $\omega=g(J\cdot, \cdot)$ in $\Lambda^2(M)$. Since $\omega$ is everywhere non-degenerate the manifold $M$ is naturally
oriented by the top degree form $\omega^m$. \\
We shall give now a short account on some of the $U(m)$-modules of relevance for our study. Let $\Lambda^{\star}$ denote the space of differential 
$p$-forms on $M$, to be assumed real-valued, unless stated otherwise. In what follows we shall use the metric to identify vectors and $1$-forms, 
almost implicitly, by means of the musical isomorphism $X \in TM \mapsto X^{\sharp} \in \Lambda^1$.

We consider the operator $\J : \Lambda^p \to \Lambda^p$ acting on 
a $p$-form $\alpha$ by
\begin{equation*}(\J \alpha)(X_1, \ldots, X_p)=\sum \limits_{k=1}^{p}\alpha(X_1, \ldots, JX_k, \ldots X_p)
\end{equation*}
for all $X_1, \ldots X_p$ in $TM$. For future use let us note that alternatively 
$$ \J \alpha=\sum \limits_{i=1}^{2n} e_i \wedge (Je_i \lrcorner \alpha)
$$ 
for all $\alpha$ in $\Lambda^{\star}(M)$, where $\{e_i, 1 \le i \le 2n\}$ is some local orthonormal frame on $M$.  The almost complex structure $J$ can also 
be extended to $\Lambda^{\star}(M)$ by setting 
$$ (J\alpha)(X_1, \ldots, X_p)=\alpha(JX_1, \ldots, JX_p)
$$
for all $X_i, 1 \le i \le p$ in $TM$. Note that $JX^{\sharp}=-(JX)^{\sharp}$ for all $X$ in $TM$, and in what follows each time the notation $JX, X \in TM$ should refer to a $1$-form it is intended to
mean $(JX)^{\sharp}$.

Now $\J$ acts as 
a derivation on $\Lambda^{\star}$ and gives the complex bi-grading of the exterior algebra in the following sense. Let $\lambda^{p,q}$ be given as the $-(p-q)^2$-eigenspace of $\J^2$. One has 
\begin{equation*}\Lambda^s=\bigoplus \limits_{p+q=s}^{} \lambda^{p,q},
\end{equation*}
an orthogonal, direct sum.  Note that $\lambda^{p,q}=\lambda^{q,p}$. Of special importance in our discussion are the spaces $\lambda^{p}=\lambda^{p,0}$;  forms $\alpha$ in $\lambda^p$ are 
such that the assignment 
$$(X_1, \ldots, X_p) \to \alpha(JX_1, X_2, \ldots, X_p)$$ 
is still an alternating form which equals $p^{-1}\J\alpha$.  We now consider the operator $L: \Lambda^{\star} \to \Lambda^{\star}$ given as multiplication with the K\"ahler form $\omega$, together with its 
adjoint, to be denoted by $L^{\star}$. Forms in 
$$ \Lambda^{\star}_0=Ker L^{\star}
$$
are called primitive and one has 
$$ [L^{\star},L]=(m-p)1_{\Lambda^p}.
$$ 
on the space of $p$-forms. Note that the inclusion $\lambda^{p} \subseteq \Lambda^{p}_0$ holds, whenever $p \ge0$.
$\\$
Let $\lambda^p \ho \lambda^q$ be the space of tensors 
$Q\in\lambda^p \otimes\lambda^q$ which satisfy 
$$Q \bJ =-\bJ Q
$$ 
where $\bJ$ as a map of $\lambda^p$ stands in fact for 
$p^{-1}\J$.  We also define $\lambda^p \ko \lambda^q$ as the space of tensors $Q : \lambda^p \to \lambda^q$ such that $Q \bJ=\bJ Q$. The following Lemma 
is easy to verify. 
\begin{lema}\label{l31} Let $a : \lambda^p \otimes \lambda^q \to \Lambda^{p+q}$ be the total alternation map. Then: 
\begin{itemize}
\item[(i)]
the image of the restriction of a to $ \lambda^p \ho \lambda^q \to \Lambda^{p+q}$ is contained in $\lambda^{p,q}$;
\item[(ii)]
the image of the restriction of a to $ \lambda^p \ko \lambda^q \to \Lambda^{p+q}$ is contained in $\lambda^{p+q}$;
\item[(iii)]
the total alternation map $a : \lambda^p \ho \lambda^q \to \Lambda^{p+q}$ is injective for $p \neq q$;
\item[(iv)]
the kernel of $a : \lambda^p \otimes \lambda^q \to \Lambda^{p+q}$ is contained in $\lambda^p \otimes_{2} \lambda^q$.
\end{itemize}
\end{lema}
In this paper most of our computations will involve forms of degree up to $4$. For further use we recall that one has 
$$ \J=J
$$
on $\lambda^{1,2}$ and also that 
\begin{equation*}
J_{\lambda^{2,2} \oplus \lambda^4}=1_{\lambda^{2,2} \oplus \lambda^4}, \ J_{\lambda^{1,3}}=-1_{\lambda^{1,3}}.
\end{equation*}
We will now briefly introduce a number of algebraic operators which play an important r\^ole in the study of connections with totally 
skew-symmetric torsion. If $\alpha$ belongs to $\Lambda^2(M)$ and 
$\varphi$ is in $\Lambda^{\star}(M)$ we define the commutator \footnote{this is proportional to the commutator in the Clifford algebra bundle $Cl(M)$ of $M$.}
\begin{equation*}
[\alpha,\varphi]=\sum \limits_{i=1}^{2n} (e_i \lrcorner \alpha) \wedge (e_i \lrcorner \varphi),
\end{equation*}
where $\{e_i, 1 \le i \le 2n \}$ is a local orthonormal frame in $TM$. Note that 
$$ [\omega,\varphi]=-\J \varphi
$$
for all $\varphi$ in $\Lambda^{\star}$. It is also useful to mention that:
\begin{equation} \label{comrule}
\begin{split}
& [\lambda^{1,1}, \lambda^3] \subseteq \lambda^3, \ [\lambda^{1,1}, \lambda^{1,2}] \subseteq \lambda^{1,2} \\
& [\lambda^2, \lambda^3] \subseteq \lambda^{1,2}. 
\end{split}
\end{equation} 
Lastly, for any couple of forms $(\varphi_1, \varphi_2)$ in $\Lambda^p \times \Lambda^q$ we define their product $\varphi_1 \b \varphi_2$ in $\Lambda^{p+q-2}$ by 
\begin{equation*}
\varphi_1 \b \varphi_2=\sum \limits_{k=1}^{2m} (e_k \lrcorner \varphi_1) \wedge (e_k \lrcorner \varphi_2),
\end{equation*}
for some local orthonormal frame $\{e_k, 1 \le k \le 2m \}$ on $M$. We end this section by listing a few properties of the product $\b$ for low degree forms, to be used in what follows.
\begin{lema} \label{btype}
The following hold:
\begin{itemize}
\item[(i)] $\varphi_1 \b \varphi_2$ belongs to $\lambda^{2,2} \oplus \lambda^{1,3}$ for all $\varphi_1, \varphi_2$ in $\lambda^{1,3}$;
\item[(ii)] $\varphi \b J\varphi$ is in $\lambda^{1,3}$ for all $\varphi$ in $\lambda^{1,2}$;
\item[(iii)] $\varphi \b \psi$ belongs to $\lambda^{1,3} \oplus \lambda^4$ whenever $\varphi$ is in $\lambda^{1,3}$ and all $\psi$ in $\lambda^3$;
\item[(iv)] $\psi_1 \b \psi_2$ belongs to $\lambda^{2,2} \oplus \lambda^4$ for all $\psi_1, \psi_2$ in $\lambda^3$;
\item[(v)] if $\psi$ is in $\lambda^3$ we have that $\psi \b J \psi=0$.
\end{itemize}
\end{lema}
The proofs consist in a simple verification which is left to the reader. 
%%%%%%%%%%%%%%%%%%%%%%%%%%%%%%%%%%%%%%%%%%%
\subsection{Curvature tensors in the almost Hermitian setting}
We shall present here a number of basic facts concerning algebraic curvature tensors we shall need later on. Since this is intended to be mainly at an algebraic 
level, our context will be that of a given Hermitian vector space $(V^{2m},g,J)$. Let us recall that the Bianchi map $b_1 : \Lambda^2 \otimes \Lambda^2 \to \Lambda^1 \otimes \Lambda^3 $ 
is defined by 
\begin{equation*}
(b_1R)_x=\sum \limits_{i=1}^{2m} e_i \wedge R(e_i, x)
\end{equation*}
for all $R$ in $\Lambda^2 \otimes \Lambda^2$ and for all $x$ in $V$. The space of algebraic curvature tensors on $V$ is given by 
\begin{equation*}
\K(\mathfrak{so}(2m))=(\Lambda^2 \otimes \Lambda^2 ) \cap Ker(b_1)
\end{equation*}
and it is worth observing that $\K(\mathfrak{so}(2m))=S^2(\Lambda^2) \cap Ker(a)$. Restricting the group to the unitary one makes appear the space of K\"ahler curvature 
tensors given by 
\begin{equation*}
\K(\mathfrak{u}(m))= (\lambda^{1,1} \otimes \lambda^{1,1}) \cap Ker(b_1)=S^2(\lambda^{1,1}) \cap Ker(a).
\end{equation*}
We shall briefly present some well known facts related to the space of K\"ahler curvature tensors, with proofs given mostly for the sake of completeness. First of all, we have an inclusion of $\lambda^{2,2}$ into 
$S^2(\lambda^{1,1})$ given by 
$$ \Omega \mapsto \hat{\Omega}
$$
where $\hat{\Omega}(x,y)=-\frac{1}{4} (\Omega(x,y)+\Omega(Jx,Jy))$ for all $x,y$ in $V$. This is mainly due to the fact that forms of type $(2,2)$ are $J$-invariant. A short computation yields 
\begin{equation} \label{b1f}
b_1(\hat{\Omega})=\Omega,
\end{equation}
for all $\Omega$ in $\lambda^{2,2}$. There is also an embedding of $\lambda^{1,3}$ into $\lambda^2 \otimes \lambda^{1,1}$ given by $\Omega \mapsto \check{\Omega}$ where 
$$ \check{\Omega}(x,y)=\Omega(x,y)-\Omega(Jx,Jy)
$$
for all $x,y$ in $V$. It is easily verified that 
\begin{equation} \label{b2ff}
(b_1 \check{\Omega})_x=Jx \lrcorner \J \Omega-2 x \lrcorner \Omega
\end{equation}
for all $x,y$ in $V$ and for all $\Omega$ in $\lambda^{1,3}$. 

Similarly to the well known splitting $S^2(\Lambda^2)=\K(\mathfrak{so}(2m)) \oplus \Lambda^4$ we 
have:
\begin{pro} \label{splitK1}
There is an orthogonal, direct sum splitting 
$$ S^2(\lambda^{1,1})=\K(\mathfrak{u}(m)) \oplus \lambda^{2,2}.
$$
Explicitly, any $Q$ in $ S^2(\lambda^{1,1})$ can be uniquely written as $Q=R+\hat{\Omega}$ for some K\"ahler curvature tensor $R$ and some $\Omega$ in $\lambda^{2,2}$. 
\end{pro}
\begin{proof}
Let $Q$ belong to $S^2(\lambda^{1,1})$. It satisfies $Q(Jx,Jy)=Q(x,y)$ for all $x,y$ in $V$ and since $Q$ belongs to $S^2(\Lambda^2)$ it splits as 
$$ Q=R+\Omega
$$
where $R$ is in $\K(\mathfrak{so}(2m))$ and $\Omega$ belongs to $\Lambda^4$. In particular 
\begin{equation} \label{eqh}
R(Jx,Jy)+\Omega(Jx,Jy)=R(x,y)+\Omega(x,y)
\end{equation}
for all $x,y$ in $V$. \\
Since $Q$ belongs to $\lambda^{1,1} \otimes \lambda^{1,1}$ we have that $(b_1Q)_x$ belongs to $\lambda^{1,2}$ for all $x$ in $V$. But $R$ is a curvature tensor, hence
$(b_1Q)_x=-3x \lrcorner \Omega$, leading to $x \lrcorner \Omega$ in $\lambda^{1,2}$ for all $x$ in $V$. It is then easy to see that $\Omega$ must be an element of $\lambda^{2,2}$. It follows by direct verification that 
$R_{\Omega}$ given by 
$$ R_{\Omega}(x,y)=\Omega(Jx,Jy)-\frac{1}{3}\Omega(x,y)
$$
belongs to $\K(\mathfrak{so}(2m))$. Setting now $R_0=R-\frac{3}{4}R_{\Omega}$ it is easy to see from \eqref{eqh} that $R_0$ belongs to $\K(\mathfrak{u}(m))$ and our claim follows after appropriate rescaling from 
$$ Q(x,y)=R_0(x,y)+\frac{3}{4}\biggl [ \Omega(Jx,Jy)+\Omega(x,y) \biggr ]
$$
for all $x,y$ in $V$. 
\end{proof}
Our next and last goal in this section is to have an explicit splitting of certain elements of $\Lambda^2 \otimes \lambda^{1,1}$ along $\Lambda^2=\lambda^{1,1} \oplus \lambda^2$. More explicitly, we will look at this 
in the special case of a tensor $\RR$ in $\Lambda^{2} \otimes \lambda^{1,1}$ such that 
\begin{equation} \label{curvvh}
\RR(x,y,z,u)-\RR(z,u,x,y)=\gamma_x(y,z,u)-\gamma_y(x,z,u)-\gamma_z(u,x,y)+\gamma_u(z,x,y)
\end{equation} 
holds for all $x,y,z,u$ in $V$. Here $\gamma$ belongs to $\Lambda^1 \otimes \Lambda^3$ and in what follows we shall use the notation 
$$ T=a(\gamma)
$$
for the total alternation of $\gamma$. The tensor $\RR$ is the algebraic model for the curvature tensor of a Hermitian connection with 
totally skew-symmetric torsion, which will be our object of study later on in the paper.
\begin{rema} \label{trivanheckepaper}
A complete decomposition of $\mathcal{K}(\mathfrak{so}(2m))$ into irreducible components under the action of $U(m)$ has been given by Tricerri and Vanhecke in \cite{Tri}.  Further information 
concerning the splitting of $\Lambda^2 \otimes \Lambda^2$, again under the action of $U(m)$ is given in detail in \cite{Falc}.  While the material we shall present next can be equivalently derived from 
these references, it is given both for self-containedness and also as an illustration that one can directly proceed, in the case of a connexion with torsion, to directly split its curvature tensor without reference to the 
Riemann one. As a general observation we also note that the procedure involves only control of the orthogonal projections onto the relevant $U(m)$-submodules. 

\end{rema}

Recall that for any $\alpha$ in $\Lambda^2$ its orthogonal projection on $\lambda^{1,1}$ is given by
$$ \alpha_{\lambda^{1,1}}=\frac{1}{2}(\alpha+J\alpha).
$$
To obtain the decomposition of $\RR$ we need the following preliminary Lemma. 
\begin{lema} \label{b1gen}
Let $\eta$ in $\lambda^1 \otimes \lambda^{1,2}$ be given. We consider the tensors 
$H_1, H_2$ in $\Lambda^2 \otimes \lambda^{1,1}$ given by 
\begin{equation*}
\begin{split}
H_1(x,y)=& (x \lrcorner \eta_y-y \lrcorner \eta_x)_{\lambda^{1,1}}\\
H_2(x,y)=& H_2(Jx,Jy)
\end{split}
\end{equation*}
for all $x,y$ in $V$. Then:
\begin{itemize}
\item[(i)] $2(b_1 H_1)_x=3\eta_x+J\eta_{Jx}+x \lrcorner a(\eta)-Jx \lrcorner a(J\eta)$; 
\item[(ii)] $2(b_1 H_2)_x=\eta_x+3 J\eta_{Jx}+Jx \lrcorner a^c(\eta)+x \lrcorner a^c(J\eta)$ 
\end{itemize}
for all $x$ in $V$. Here $J\eta$ in $\lambda^1 \otimes \lambda^{1,2}$ is defined by $(J\eta)_x=J\eta_x$ for all $x$ in $V$, and the complex alternation 
map $a^c : \lambda^1 \otimes \lambda^{1,2} \to \Lambda^4$ is given by 
\begin{equation*}
a^c(\z)=\sum \limits_{k=1}^{2m} e_k \wedge \eta_{Je_k}
\end{equation*}
\end{lema}
\begin{proof}
(i) Let $\{e_k, 1 \le k \le 2m\}$ be an orthonormal basis in $V$. Directly from its definition, the tensor $H_1$ is given by 
$$ 2H_1(x,y)=x \lrcorner \eta_{y}-Jx \lrcorner J\eta_{y}-y \lrcorner \eta_{x}+Jy \lrcorner J\eta_{x}
$$
hence 
$$ 2 (b_1H_1)_x=4\eta_x-\sum \limits_{k=1}^{2m} e_k \wedge (x \lrcorner \eta_{e_k})+\sum \limits_{k=1}^{2m} e_k \wedge (Jx \lrcorner J\eta_{e_k})
$$
since $\J=J$ on $\lambda^{1,2}$. Now 
$$ x \lrcorner a(\eta)=\eta_x-\sum \limits_{k=1}^{2m} e_k \wedge (x \lrcorner \eta_{e_k})$$
and 
$$ Jx \lrcorner a(J\eta)=J\eta_{Jx}-\sum \limits_{k=1}^{2m} e_k \wedge (Jx \lrcorner J\eta_{e_k}) $$
for all $x$ in $V$ and the claim follows immediately.\\
(ii) it is enough to use (i) when replacing $\eta$ by $J\eta$.
\end{proof}
We split now 
$$ \gamma=\gamma^{1,2}+\gamma^3
$$
along $\Lambda^1 \otimes \Lambda^3=(\Lambda^1 \otimes \lambda^{1,2}) \oplus (\Lambda^1 \otimes \lambda^3)$ and also 
$$T=T^{2,2}+T^{1,3}+T^4$$ 
along the bidegree decomposition of $\Lambda^4$ in order to be able to introduce some of the components of the tensor $\RR$.
These are the tensor $\RR^a$ in $\lambda^{1,1} \otimes \lambda^{1,1}$ given by
\begin{equation*}
\RR^a(x,y)=(\gamma_xy-\gamma_yx)_{\lambda^{1,1}}+(\gamma_{Jx}Jy-\gamma_{Jy}Jx)_{\lambda^{1,1}}-\frac{1}{2}(T^{2,2}(x,y)+T^{2,2}(Jx,Jy)),
\end{equation*}
for all $x,y$ in $V$, and the tensor $\RR^m$ in $\lambda^2 \otimes \lambda^{1,1}$ defined by 
\begin{equation*}
\RR^m(x,y)=(\gamma_xy-\gamma_yx)_{\lambda^{1,1}}-(\gamma_{Jx}Jy-\gamma_{Jy}Jx)_{\lambda^{1,1}}-\frac{1}{2}(T^{1,3}(x,y)-T^{1,3}(Jx,Jy)),
\end{equation*}
for all $x,y$ in $V$. Here we use $\gamma_xy$ as a shorthand for $y \lrcorner \gamma_x$, whenever $x,y$ are in $V$. 

The promised decomposition result for $\RR$ is achieved mainly by projection of \eqref{skewpart} onto 
$$ \Lambda^2 \otimes \lambda^{1,1}=(\lambda^{1,1} \otimes \lambda^{1,1}) \oplus (\lambda^{1,1} \otimes \lambda^2)
$$
while taking into account that $\RR$ belongs to $\Lambda^2 \otimes \lambda^{1,1}$. 
\begin{teo} \label{curvfin0}Let $\RR$ belong to $\Lambda^2 \otimes \lambda^{1,1}$ such that \eqref{curvvh} is satisfied. We have a decomposition 
\begin{equation*}
\RR=\RR^K+\hat{\Omega}+\frac{1}{2}\RR^a+\RR^m
\end{equation*}
where $\RR^K$ belongs to $\mathcal{K}(\mathfrak{u}(m))$ and $\Omega$ is in $\lambda^{2,2}$ . Moreover:
\begin{itemize}
\item[(i)] $R^a$ belongs to $\Lambda^{2}(\lambda^{1,1})$ and satisfies the Bianchi identity  
\begin{equation*}
(b_1 R^a)_x=-2(\gamma^{1,2}_x+J\gamma_{Jx}^{1,2})-\frac{1}{2}x \lrcorner A_1-\frac{1}{2} Jx \lrcorner A_2,\\
\end{equation*}
for all $x$ in $V$. The forms $A_1$ and $A_2$ are explicitly given by 
\begin{equation*}
\begin{split}
A_1=&a(\gamma^{1,2})+a^c(J\gamma^{1,2})-4T^{2,2},\\
A_2=&a^c(\gamma^{1,2})-a(J\gamma^{1,2}).
\end{split}
\end{equation*}
\item[(ii)] the Bianchi identity for $R^m$ in $\lambda^2 \otimes \lambda^{1,1}$ is 
\begin{equation*}
(b_1 R^m)_x=-(\gamma^{1,2}_x-J\gamma^{1,2}_{Jx})-\frac{1}{2}x \lrcorner A_3+\frac{1}{2}Jx \lrcorner A_4,\\
\end{equation*}
for all $x$ in $V$, where 
\begin{equation*}
\begin{split}
A_3=&a(\gamma^{1,2})-a^c(J\gamma^{1,2})-2T^{1,3},\\
A_4=&a(J\gamma^{1,2})+a^c(\gamma^{1,2})-\J T^{1,3}.
\end{split}
\end{equation*}
\end{itemize}
\end{teo}
\begin{proof}
We first show that $R^a$ belongs to $\Lambda^2(\lambda^{1,1})$. 
Directly from the definition of the map $a$ we have that the tensor 
$$ Q(x,y):=\gamma_xy-\gamma_yx-\frac{1}{2}T(x,y)
$$
belongs to $\Lambda^2(\Lambda^2)$, so after projection on $\lambda^{1,1} \otimes \lambda^{1,1}$ it follows that 
$$ (\gamma_xy-\gamma_yx)_{\lambda^{1,1}}+(\gamma_{Jx}Jy-\gamma_{Jy}Jx)_{\lambda^{1,1}}-\frac{1}{2}(T(x,y)+T(Jx,Jy))_{\lambda^{1,1}}$$
is an element of $\Lambda^2(\lambda^{1,1})$. We conclude by recording that 
$$ (T(x,y)+T(Jx,Jy))_{\lambda^{1,1}}=T^{2,2}(x,y)+T^{2,2}(Jx,Jy)
$$
since forms in $\lambda^{1,3}$ are $J$-anti-invariant whilst those in $\lambda^{2,2} \oplus \lambda^4$ are $J$-invariant and moreover any $\varphi$ in 
$\lambda^4$ satisfies $\varphi(J \cdot, J \cdot, \cdot, \cdot)=-\varphi(\cdot, \cdot, \cdot, \cdot)$. \\
The next step is to notice that \eqref{curvvh} actually says that $\RR-Q$ belongs to $S^2(\Lambda^2)$ and to project again on $\lambda^{1,1} \otimes \lambda^{1,1}$. By the argument above 
it follows that 
$$ \RR(x,y)+\RR(Jx,Jy)-R^a(x,y)$$
is in $S^2(\lambda^{1,1})$ thus by Proposition \ref{splitK1} we can write 
\begin{equation} \label{step1}
\RR(x,y)+\RR(Jx,Jy)-R^a(x,y)=2R^K(x,y)+2\hat{\Omega}(x,y),
\end{equation}
for all $x,y$ in $V$, where $\Omega$ is in $\lambda^{2,2}$ and $R^K$ belongs to $\K(\mathfrak{u}(m))$.\\
Now, rewriting \eqref{curvvh} as 
\begin{equation*}
\RR(x,y,z,u)-\RR(z,u,x,y)=2(\gamma_x(y,z,u)-\gamma_y(x,z,u))-T(x,y,z,u)
\end{equation*}
we obtain after projection on $\lambda^2 \otimes \lambda^{1,1}$ that 
\begin{equation*}
\RR(x,y)-\RR(Jx,Jy)=2(\gamma_xy-\gamma_yx)_{\lambda^{1,1}}-2(\gamma_{Jx}Jy-\gamma_{Jy}Jx)_{\lambda^{1,1}}-(T(x,y)-T(Jx,Jy))_{\lambda^{1,1}}
\end{equation*}
for all $x,y$ in $V$. Since $\lambda^4 \oplus \lambda^{2,2}$ consists in $J$-invariant forms, it is easy to see that 
$$
(T(x,y)-T(Jx,Jy))_{\lambda^{1,1}}=T^{1,3}(x,y)-T^{1,3}(Jx,Jy),
$$
in other words
\begin{equation} \label{step2} 
\frac{1}{2}(\RR(x,y)-\RR(Jx,Jy))=\RR^m(x,y)
\end{equation}
for all $x,y$ in $V$. The splitting of $\RR$ follows now from \eqref{step1} and \eqref{step2} and from the fact that the component of $\gamma$ on $\lambda^1 \otimes \lambda^3$ is not seen 
by the projection on $\lambda^{1,1}$. 

To finish the proof, it remains only to prove the Bianchi identities for $\RR^a$ and $\RR^m$. Using Lemma \ref{b1gen} and \eqref{b1f} it is easy to get that 
\begin{equation*}
\begin{split}
(b_1 R^a)_x=&-2(\gamma^{1,2}_x+J\gamma_{Jx}^{1,2})+2x \lrcorner T^{2,2}\\
&-\frac{1}{2}x \lrcorner (a(\gamma^{1,2})+a^c(J\gamma^{1,2}))-\frac{1}{2}Jx \lrcorner (a^c(\gamma^{1,2})-a(J\gamma^{1,2}))
\end{split}
\end{equation*}
for all $x$ in $V$. The Bianchi identity for $\RR^m$ is proved along the same lines and therefore left to the reader.
\end{proof}
\begin{rema} 
Underlying Theorem \ref{curvfin0} are the following isomorphisms of $\mathfrak{u}(m)$-modules. The first is 
$$ b_1 : \Lambda^2(\lambda^{1,1}) \to (\lambda^1 \otimes_1 \lambda^{1,2} ) \cap Ker(a),
$$
where 
$$ \lambda^1 \otimes_1 \lambda^{1,2}=\{\z \in \lambda^1 \otimes \lambda^{1,2}: \z_{Jx}=-J \z_x \ \mbox{for all} \ x \in V\}.
$$
The second is given by 
$$ b_1 : \lambda^2 \otimes \lambda^{1,1} \to \lambda^1 \otimes_2 \lambda^{1,2} 
$$
where $ \lambda^1 \otimes_2 \lambda^{1,2}=\{\z \in \lambda^1 \otimes \lambda^{1,2}: \z_{Jx}=J \z_x \ \mbox{for all} \ x \in V\}$.\\
This makes that in practice it is not necessary to work with the somewhat involved expressions for $R^a,R ^m$ but rather with their Bianchi contractions, which 
are tractable.
\end{rema}
%%%%%%%%%%%%%%%%%%%%%%%%%%%%%%%%%%%%%%%%%%%%%%%%%%%%%%%%%%%%%%%%%%%%%
\subsection{The Nijenhuis tensor and Hermitian connexions}
Let $(M^{2m},g,J)$ be an almost Hermitian manifold. Recall that the Nijenhuis tensor of the almost complex structure $J$ is defined by 
$$ N_J(X,Y)=[X,Y]-[JX,JY]+J[JX,Y]+J[X,JY]
$$
for all vector fields $X,Y$ on $M$. It satisfies 
\begin{equation*}
\begin{split}
&N_J(X,Y)+N_J(Y,X)=0\\
&N_J(JX,JY)=-N_J(X,Y)\\
&N_J(JX,Y)=-JN_J(X,Y)
\end{split}
\end{equation*}
for all $X,Y$ in $TM$. By evaluating $N_J$ using the Riemannian metric $g$ we can also form the tensor $N^J$ defined by 
$$ N^J_X(Y,Z)=g(N_J(Y,Z),X)
$$
for all $X,Y,Z$ in $TM$, which therefore belongs to $\lambda^1 \otimes_2 \lambda^2$.  
When $N_J$ vanishes identically, $J$ is said to be integrable and gives $M$ the structure of a complex manifold. \\
Denoting by $\nabla$ the Levi-Civita connexion attached to the Riemannian 
metric $g$ we form the tensor $\nabla J$ in $\lambda^1 \otimes \lambda^2$. Then one can alternatively compute the Nijenhuis tensor as 
\begin{equation} \label{nlc}
N_J(X,Y)=-\biggl [ (\nabla_{JX}J)Y-(\nabla_{JY}J)X\biggr ]+J\biggl [ (\nabla_XJ)Y- (\nabla_YJ)X\biggr ].
\end{equation}
It is now a good moment to recall that an almost Hermitian structure $(g,J)$ is called K\"ahler if and only if $\nabla J=0$, or equivalently $\nabla \omega=0$. This implies the integrability of $J$ by making use of \eqref{nlc} and  also that 
$\omega$ is a symplectic form, in the sense that $d \omega=0$ where $d$ denotes the exterior derivative. It turns out that for an arbitrary almost complex structure $(g,J)$, the tensors 
$d \omega$ in $\Lambda^3(M)$ and $N^J$ in $\lambda^1 \otimes \lambda^2$ form a full set of obstructions to having $(g,J)$ K\"ahler as the following general fact shows.
\begin{pro} \label{st1}
For any almost Hermitian structure $(g,J)$ on $M$ we have
\begin{equation} \label{jet1}
2\nabla_X \omega=-N^J_{JX}+X \lrcorner d \omega+JX \lrcorner Jd \omega,
\end{equation}
for all $X$ in $TM$. 
\end{pro}
\begin{proof}Although this is a standard fact we give the proof for self-containedness. From the definition of the exterior derivative we have 
\begin{equation*}
d \omega (X,Y,Z)=(\nabla_X \omega)(Y,Z)-(\nabla_Y \omega)(X,Z)+(\nabla_Z \omega)(X,Y)
\end{equation*}
for all $X,Y,Z$ in $TM$. Since $\nabla_X \omega$ is in $\lambda^2$ for all $X$ in $TM$ we obtain 
\begin{equation*}
\begin{split}
&d\omega(X,Y,Z)-d\omega(JX,JY,Z)\\
&=(\nabla_X \omega)(Y,Z)-(\nabla_Y \omega)(X,Z)-(\nabla_{JX} \omega)(JY,Z)+(\nabla_{JY} \omega)(JX,Z)\\
& +2(\nabla_Z \omega)(X,Y)
\end{split}
\end{equation*}
for all $X,Y,Z$ in $TM$. This can be rewritten by using \eqref{nlc} as 
\begin{equation*}
d\omega(X,Y,Z)-d\omega(JX,JY,Z)=\langle N_J(X,Y), JZ\rangle +2(\nabla_Z \omega)(X,Y)
\end{equation*}
whenever $X,Y,Z$ belong to $TM$, and the claim follows.
\end{proof}
Thus $d\omega=0$ and $N_J=0$ result in $\nabla J=0$, in other words in $(g,J)$ being K\"ahler.\\
We shall call a linear connection on the tangent bundle to $M$ Hermitian if it respects both the metric and the almost complex structure. In the framework of almost Hermitian geometry two connections 
play a distinguished r\^ole. The first is called the first canonical Hermitian connection and it is defined by 
$$ \bnabla_X=\nabla_X+\eta_X
$$
where $\eta$ in $\lambda^1 \otimes \lambda^2$ is given by $\eta_X=\frac{1}{2}(\nabla_X J)J$. The tensor $\eta$ is called the intrinsic torsion tensor of the $U(m)$-structure on $M$ induced by 
$(g,J)$. The connexion $\bnabla$ is minimal \cite{pg}, in the sense of minimising the norm within the space of almost Hermitian connections. \\
To introduce the second Hermitian connexion we need some preliminaries. First of all let us decompose the $3$-form $d\omega=d^{1,2}\omega+d^3\omega$ along the splitting 
$\Lambda^3(M)=\lambda^{1,2} \oplus \lambda^3$. Then:
\begin{lema} \label{anij} 
Let $(M^{2m},g,J)$ be almost Hermitian. Then $a(N^J)=4Jd^3\omega$.
\end{lema}
\begin{proof} Using \eqref{jet1} we obtain
$$ d\omega=a(\nabla \omega)=-\sum \limits_{i=1}^{2m}e_i \wedge N^J_{Je_i}+3d\omega+\J(J d\omega)
$$
where $\{e_i, 1 \le i \le 2n\}$ is some local orthonormal frame. Since $N^J$ belongs to $\lambda^1 \otimes_2 \lambda^2$ it is easy to see that 
$ \J(a(N^J))=3\sum \limits_{i=1}^{2m}e_i \wedge N^J_{Je_i}$ and the claim follows by taking into account that $\J=J$ on $\lambda^{1,2}$.
\end{proof}
Since $\lambda^3 \subseteq \lambda^1 \otimes_2 \lambda^3$ it follows that the Nijenhuis tensor splits as 
\begin{equation*}
N^J_X=\hat{N}^J_X+\frac{4}{3}X \lrcorner Jd^3\omega
\end{equation*}
for all $X$ in $TM$, where the tensor $\hat{N}^J$ belongs to the irreducible $U(m)$-module 
$$W_2=(\lambda^1 \otimes_2 \lambda^2) \cap ker(a).$$
\begin{pro} \label{c2}
Let $(M^{2m},g,J)$ be almost Hermitian. The linear connexion $D$ defined by $D_X=\nabla_X+\z_X$ where 
$$ 2\z_X=X \lrcorner Jd^{1,2} \omega-\frac{1}{3} X \lrcorner Jd^3\omega+\frac{1}{2}\hat{N}^J_X
$$
is almost Hermitian. 
\end{pro}
\begin{proof} That $D$ is metric is clear since $\z_X$ is a two form for all $X$ in $TM$. To see that $D\omega=0$ it is enough to 
show that $\nabla_X\omega+[\z_X,\omega]=0$ or equivalently 
$$ \nabla_X\omega+\J\z_X=0
$$ 
for all $X$ in $TM$. But 
\begin{equation*}
\begin{split}
2\J \z_X=&\J( X \lrcorner Jd^{1,2}\omega)-\frac{1}{3} \J(X \lrcorner Jd^3\omega)+\frac{1}{2} \J \hat{N}^J_X\\
=&-X \lrcorner d^{1,2} \omega-JX \lrcorner Jd^{1,2}\omega-\frac{2}{3}X \lrcorner d^3 \omega+\hat{N}^J_{JX}
\end{split}
\end{equation*}
after taking into account that $\hat{N}^J$ belongs to $\lambda^1 \otimes_2 \lambda^2$. Therefore, after taking into account \eqref{jet1}
\begin{equation*}
\begin{split}
2\nabla_X\omega+2\J \z_X=&-\biggl [ \hat{N}^J_{JX}+ \frac{4}{3}JX \lrcorner Jd^3\omega \biggr ]+X \lrcorner d \omega+JX \lrcorner Jd \omega\\
&-X \lrcorner d^{1,2} \omega-JX \lrcorner Jd^{1,2}\omega-\frac{2}{3}X \lrcorner d^3 \omega+\hat{N}^J_{JX}=0
\end{split}
\end{equation*}
after a straightforward computation, and the result follows.
\end{proof}
It follows from Proposition \ref{st1} and Lemma \ref{anij} that the intrinsic torsion tensor of the almost Hermitian structure $(g,J)$ is completely 
determined by 
$$ (d\omega, \hat{N}^J) \ \mbox{in} \ \Lambda^3 \oplus W_2.
$$
Since $\Lambda^3 \oplus W_2$ has four irreducible components under the action of $U(m)$ the Gray-Hervella classification \cite{GrayH} singles out $16$-classes of almost Hermitian manifolds. Similar classification 
results are available for the groups $G_2$ \cite{FGray, CsIv2} and $Spin(7)$ \cite{Fern} as well as for quaternionic structures \cite{CaSw} and 
$SU(3)$-structures on $6$-dimensional manifolds \cite{ChSa}. The case of $Spin(9)$-structures on $16$-dimensional manifolds 
has been equally treated in \cite{Fri9}. For $PSU(3)$-structures on $8$-dimensional manifolds and $SO(3)$-structures in dimension $5$ the decomposition of the intrinsic torsion tensor has been studied 
in \cite{Hit1, Witt1} and \cite{Pawel}.
%%%%%%%%%%%%%%%%%%%%%%%%%%%%
\subsection{Admissible totally skew-symmetric torsion}
In this section we shall start to specialise our discussion to a particular class of almost Hermitian manifolds, to be characterised 
in terms of the torsion tensor of the Hermitian connexion $D$. We recall that the torsion tensor of a linear connexion on $M$, e.g. $D$, is the tensor $T^D$ in $\Lambda^2 \otimes \Lambda^1$ given by 
$$ T^D(X,Y)=D_XY-D_YX-[X,Y]
$$
for all vector fields $X,Y$ on $M$. In the case of the connexion $\bnabla$, the torsion will be denoted simply by $T$. \\
The following result of Friedrich and Ivanov clarifies 
in which circumstances an almost Hermitian structure admits a Hermitian connexion with totally skew-symmetric torsion.
\begin{teo} \label{caract} \cite{FriIva} Let $(M^{2m},g,J)$ be almost Hermitian. There exists an almost Hermitian connexion with torsion in $\Lambda^3(M)$ if and 
only if $N^J$ is a $3$-form. If the latter holds, the connexion is unique and equals $D$.
\end{teo} 
Almost Hermitian manifolds with totally skew-symmetric Nijenhuis tensor form the so-called Gray-Hervella class $W_1+W_3+W_4$ and are usually 
called $\mathcal{G}_1$ almost Hermitian structures. As it follows from the discussion above an almost Hermitian structure $(g,J)$ belongs to 
the class $\mathcal{G}_1$ if and only if $\hat{N}^J=0$. Because of its unicity, the connexion $D$ in Theorem \ref{caract} will be referred to as the characteristic connexion of the $\mathcal{G}_1$-manifold 
$(M^{2m},g,J)$. 

In the remaining of this section we will work on a given $\mathcal{G}_1$-manifold $(M^{2m},g,J)$ and we will derive a number of facts to be used further on. 
We have: 
\begin{equation*}
\begin{split}
&N^J=\frac{4}{3}Jd^3\omega\\
&2\z=T^D=Jd^{1,2}\omega-\frac{1}{3}Jd^3\omega.
\end{split}
\end{equation*}
It is also worth noting that \eqref{jet1} is then updated to:
\begin{equation*}
2\nabla_X\omega=X \lrcorner d^{1,2} \omega+JX \lrcorner Jd^{1,2} \omega+\frac{2}{3} X \lrcorner d^3\omega
\end{equation*}
for all $X$ in $TM$. Under a shorter and perhaps more tractable form this reads
\begin{equation} \label{g1}
\nabla_X\omega=X \lrcorner t+JX \lrcorner Jt+X \lrcorner \psi^{+}
\end{equation}
for all $X$ in $TM$, where the $3$-forms $t$ in $\lambda^{1,2}$ and $\psi^{+}$ in $\lambda^3$ are given by 
$$ t=\frac{1}{2}d^{1,2}\omega , \ \psi^{+}=\frac{1}{3}d^3\omega .$$
In the same spirit we can also re-express $d\omega$ and the torsion form as 
\begin{equation} \label{do}
d\omega=2t+3\psi^{+}
\end{equation}
and 
\begin{equation} \label{tD}
T^D=2Jt+\psi^{-}
\end{equation}
where the $3$-form $\psi^{-}$ in $\lambda^3$ is given by $\psi^{-}=\psi^{+}(J \cdot, \cdot, \cdot)$.
\begin{rema} (i)When $m=2$ any almost Hermitian structure of class $\mathcal{G}_1$ is automatically Hermitian, that is $N_J=0$. This is due to the vanishing 
of $\lambda^3$ in dimension $4$. In what follows we shall therefore assume that $m \ge 3$.\\
(ii) From \eqref{g1} it easy to see that almost Hermitian structures in the class $\mathcal{G}_1$ are alternatively described as those satisfying $(\nabla_{JX}J)JX=-(\nabla_XJ)X$ for all 
$X$ in $TM$.
\end{rema}
When $J$ is integrable it is easy to see that $T^D=2Jd\omega$ belongs to $\Lambda^{3}(M)$ (actually $\lambda^{1,2}$ after taking into account Lemma \ref{anij}). 
In this case the connection $D$ is referred to as the Bismut connexion \cite{bismut}.

Almost quaternionic or almost hyperhermitian structures admitting a structure preserving connexion with totally skew-symmetric torsion have been introduced and 
given various characterisations in \cite{HowePa, GraPoon, IvaQKT, CabSwann}. The resulting geometries are known under the names QKT (quaternion-K\"ahler with torsion) and HKT (hyperk\"ahler with torsion). 

%%%%%%%%%%%%%%%%%%%%%%%%%%%%%%%%%%%%%%
\subsection{Hermitian Killing forms}
In this section, briefly recalling the definition of Killing forms in Riemannian geometry we shall present a variation of that notion, more suitable 
to the almost Hermitian setting. This is done to prepare the ground to present, after proving the results in the next section, 
the relationship binding together almost Hermitian structure of type $\mathcal{G}_1$ and this kind of generalised Killing form. \\
Let $(M^{2m},g,J)$ be almost Hermitian in the class $\mathcal{G}_1$. 
To the connexion $D$ we can associate the exterior 
differential $d_D$, given as in the case of the usual exterior derivative $d$ by 
$$ d_D=\sum \limits_{i=1}^{2m} e_i \wedge D_{e_i}.
$$
Given that $D$ is a Hermitian connexion it is easy to see that 
\begin{equation} \label{dtype}
[\J,d_D]=(-1)^p Jd_DJ 
\end{equation}
holds on $\Lambda^p$. By duality, we also have that
\begin{equation} \label{dtypestar}
[\J,d_D^{\star}]=(-1)^p Jd_D^{\star}J 
\end{equation}
on $\Lambda^p$, where $d_D^{\star} : \Lambda^{\star} \to \Lambda^{\star}$ is the formal adjoint of $d_D$. Note that $d_D^{\star}$ is computed at a point $m$ of $M$ by 
$$ d_D^{\star}=-\sum \limits_{i=1}^{2m} e_i \lrcorner D_{e_i}
$$
where $\{e_i, 1 \le i\le 2m \}$ is a local frame around $m$, geodesic at $m$ w.r.t. the connection $D$.
A straightforward implication of \eqref{dtype} is that  
$$ d_D \varphi \in \lambda^{p+1,q} \oplus \lambda^{p,q+1}
$$
for all $\varphi$ in $\lambda^{p,q}$. Denoting by $\varphi_{\lambda^{p,q}}$ the orthogonal projection of $\varphi$ in $\Lambda^{\star}$ on $\lambda^{p,q}$ we split
$$ d_D=\pf_D+\bpf_D
$$
where $\pf_D \varphi=(d_D \varphi)_{\lambda^{p+1,q}}$ and $\bpf_D \varphi=(d_D \varphi)_{\lambda^{p,q+1}}$ for all $\varphi$ in $\lambda^{p,q}$.\\
%%%%%%%%%%%%%%%%%%5
%\edz{ambiguous}. \\
Lastly, we mention that 
the K\"ahler identities 
\begin{equation} \label{kid}
\begin{split}
[d_D, L^{\star}]=&(-1)^p Jd_D^{\star}J\\
[d_D^{\star}, L]=&(-1)^{p+1} Jd_DJ
\end{split}
\end{equation}
hold, whenever $\varphi$ belongs to $\Lambda^p$. We can now make the following.
\begin{defi} \label{hk}
A form $\varphi$ in $\lambda^{p,q}$ is a Hermitian Killing form if and only if 
$$ D_X \varphi=(X \lrcorner A)_{\lambda^{p,q}}
$$
for some form $A$ in $\lambda^{p,q+1} \oplus \lambda^{p+1,q}$ and for all $X$ in $TM$. 
\end{defi} 
This is much in analogy with the concept of twistor form from Riemannian geometry, see \cite{uwe1}, for details. We shall just recall that 
that a differential form $\varphi$ in $\Lambda^p(N)$, where $(N^n,h)$ is some Riemannian manifold, is called a twistor form or conformal Killing form if 
\begin{equation*}
\nabla_X\varphi= \frac{1}{p+1}X \lrcorner d\varphi+\frac{1}{n-p+1}X \wedge d^{\star} \varphi
\end{equation*}
holds, for all $X$ in $TM$. Moreover $\varphi$ is called a Killing form if it also coclosed, that is $d^{\star} \varphi=0$.
\begin{rema} Hermitian Killing $1$-forms are dual to Killing vector fields for the metric $g$. This is essentially due to the fact that that $D$ has totally skew-symmetric torsion.
\end{rema}
In the rest of this section we shall make a number of elementary observations on Hermitian Killing forms in $\lambda^{\star}$, the most 
relevant case for our aims. 
\begin{pro} \label{hk1}Let $\varphi$ in $\lambda^{p}, p \ge 2$ be a Hermitian Killing form. The following hold:
\begin{itemize}
\item[(i)] 
$ D_X \varphi=(X \lrcorner A)_{\lambda^p}$ for all $X$ in $TM$, where $A=\pf \varphi+\frac{1}{p+1}\bpf \varphi$;
\item[(ii)] $\J \varphi$ is a Hermitian Killing form.
\end{itemize}
\end{pro}
\begin{proof}
(i) From the definition  we have 
$$ D_X \varphi=(X \lrcorner A)_{\lambda^{p}}
$$
for all $X$ in $TM$, where $A$ is in $\lambda^{p,1} \oplus \lambda^{p+1}$. We split $A=B+C$ where $B$ and $C$ belong to 
$\lambda^{p,1}$ and $\lambda^{p+1}$ respectively. Since for all $X$ in $TM$:
\begin{equation*}
\begin{split}
(X \lrcorner B)_{\lambda^{p}}=&\frac{1}{2(p-1)}(JX \lrcorner \J B+(p-1) X \lrcorner B)\\
(X \lrcorner C)_{\lambda^{p}}=&X \lrcorner C
\end{split}
\end{equation*}
we find, after taking the appropriate contractions that 
\begin{equation*}
\begin{split}
d_D \varphi=&B+(p+1)C
\end{split}
\end{equation*}
and the claim follows. \\
(ii) After applying $\J$ to the Hermitian Killing equation satisfied by $\varphi$, it is enough to notice that 
$$ \J(X \lrcorner B)_{\lambda^p}=\frac{p}{p-1} (X \lrcorner \J B)_{\lambda^p}
$$
holds, together with
$$ \J (X \lrcorner C)=\frac{p}{p+1} X \lrcorner \J C.
$$
for all $X$ in $TM$.
\end{proof}
Note that, unlike Riemannian Killing forms, Hermitian Killing ones need not be necessarily coclosed. This can be easily verified by means of \eqref{dtypestar} and \eqref{kid}. The following identity is useful 
in order to better understand the notion of Hermitian Killing form. 
\begin{lema} \label{id0} Let $\varphi$ in $\lambda^p, p \ge 1$ be given. Then:
$$ D_{JX}\varphi-D_X \mathbb{J} \varphi=-2(X \lrcorner \pf_D \mathbb{J} \varphi)_{\lambda^p}$$ for all $X$ in $TM$.
\end{lema}
\begin{proof}Let us consider the tensor $\gamma$ in $\lambda^1 \otimes_1 \lambda^p$ given by 
$$ \gamma_X=D_{JX}\varphi-D_X \mathbb{J} \varphi
$$
for all $X$ in $TM$. For notational convenience we define $B=\pf_D \mathbb{J} \varphi$ in $\lambda^{1,p}$ and note, as in the proof of Proposition \ref{hk1} that 
$ X \in TM \mapsto (X \lrcorner B)_{\lambda^p}$ belongs to $\lambda^1 \otimes_1 \lambda^p$ as well. Now a straightforward calculation shows that 
$$ a(X \in TM \mapsto (X \lrcorner B)_{\lambda^p})=B.$$
At the same time one has 
\begin{equation*}
a(\gamma)=-2\pf_D \mathbb{J} \varphi
\end{equation*}
hence $\gamma_X=-2(X \lrcorner \pf_D \mathbb{J} \varphi)_{\lambda^p}$ for all $X$ in $TM$ by making use of Lemma \ref{l31}, (iii) when $p \neq 1$. When $p=1$ the claim follows by a simple 
direct verification which is left to the reader.
\end{proof}
Let us now give an equivalent  characterisation 
of Hermitian Killing forms.
\begin{pro} \label{middle}
The following hold:
\begin{itemize}
\item[(i)] a form $\varphi$ in $\lambda^p, p \ge 2$ is a Hermitian Killing form if and only if the component of $D\varphi$ on $\lambda^1 \otimes_2 \lambda^p$ is determined by $\bpf_D \mathbb{J} \varphi$, that is 
$$ D_{JX}\varphi+D_X \mathbb{J} \varphi=\frac{2}{p+1} X \lrcorner \bpf_D \mathbb{J} \varphi 
$$
for all $X$ in $TM$;
\item[(ii)] a form $\varphi$ in $\lambda^m$ is a Hermitian Killing form if and only if 
\begin{equation*}
\vert \varphi \vert D\varphi=\frac{1}{2}d \vert \varphi \vert^2 \otimes \varphi+ \frac{1}{2}Jd \vert \varphi \vert^2\otimes \mathbb{J} \varphi
\end{equation*}
\end{itemize}
\end{pro}
\begin{proof}
(i) We have 
\begin{equation*}
\begin{split}
D_X \mathbb{J} \varphi &=\frac{1}{2} (D_{JX}\varphi+D_X \mathbb{J} \varphi)-\frac{1}{2} (D_{JX}\varphi-D_X \mathbb{J} \varphi)\\
&=(X \lrcorner \pf_D \varphi)_{\lambda^3}+\frac{1}{2} (D_{JX}\varphi+D_X \mathbb{J} \varphi)
\end{split}
\end{equation*}
for all $X$ in $TM$, and the claim follows immediately from (i) in Proposition \ref{hk1}.\\
(ii) if $\varphi$ is in $\lambda^m$, we have that $\bpf_D(\mathbb{J} \varphi)=0$ since $\lambda^{m+1}=0$. From (i) we know that $\varphi$ is a Hermitian Killing form if and only if 
\begin{equation} \label{midinter}
D_{JX}\varphi+D_X \mathbb{J} \varphi=0
\end{equation}
holds for any $X$ in $TM$. Let $U$ be the open subset of $M$ given by $U=\{m \in M : \varphi_m \neq 0\}$ and let $Z$ be the complement of $U$ in $M$. We claim that \eqref{midinter} holds 
iff it holds on $U$. Indeed \eqref{midinter} holds trivially on $int(Z)$ and $U \cup int(Z)$ is dense in $M$.

Suppose now that $\varphi$ in $\lambda^m$ is a Hermitian Killing form. Then $\varphi, \mathbb{J} \varphi$ is a basis of $\lambda^m_{\vert U}$ hence 
$$ D \varphi=a \otimes \varphi+b \otimes \mathbb{J} \varphi
$$
for a couple of $1$-forms $a,b$ on $U$. The first is determined by $\vert \varphi \vert a=\frac{1}{2}d \vert \varphi \vert^2$. On the other hand the fact that $\varphi$ is a 
Hermitian form implies that $b=Ja$, hence our claim is proved on $U$ and by the density argument above on $M$. The converse statement is also clear from the previous observations.
\end{proof}
An immediate consequence of (ii) in the Proposition above is that a Hermitian Killing form in $\lambda^m$ is {\it{parallel}} with respect to the connexion $D$ as soon as it has constant length. It is now a good moment to 
provide some examples of Hermitian Killing forms.
\begin{pro} \label{example}
Let $(M^{2m},g,J)$ be a K\"ahler manifold and let $\z_k, k=1,2$ be holomorphic Killing vector fields, that is 
$$ L_{\z_k}g=0, \ L_{\z_k}J=0 \ \mbox{for} \ k=1,2.
$$
The form $\varphi=(\z_1 \wedge \z_2)_{\lambda^2}$ is a Hermitian Killing form. 
\end{pro}
\begin{proof}
The Killing equation yields 
$$ \nabla_X \z_k=\frac{1}{2}X \lrcorner d \z_k
$$
for all $X$ in $TM$ and $k=1,2$. After a few manipulations we get 
\begin{equation*}
\begin{split}
\nabla_X(\z_1 \wedge \z_2)=&\frac{1}{2}(X \lrcorner d \z_1) \wedge \z_2+\frac{1}{2} \z_1 \wedge (X \lrcorner d\z_2)\\
=& \frac{1}{2} X \lrcorner d (\z_1 \wedge \z_2)-\frac{1}{2} \langle \z_2,X\rangle d\z_1+\frac{1}{2} \langle \z_1, X \rangle d \z_2   
\end{split}
\end{equation*}
for all $X$ in $TM$. Since our given vector fields are also holomorphic, the forms $d\z_k, k=1,2$ belong to $\lambda^{1,1}$ and the claim is proved by projecting onto $\lambda^2$ while 
using the Definition \ref{hk}. We note that $\bpf \varphi=0$ because one has, as well known, $d(J\z_k)=0$ for  $k=1,2$. 
\end{proof}
We continue by giving an exterior algebra characterisation of Hermitian Killing forms, in the context of a $\mathcal{G}_1$-manifold $(M^{2m},g,J)$. 
\begin{pro} \label{id1}Let $p \ge 3$ be odd. Then:
\begin{itemize}
\item[(i)] if $\varphi$ in $\lambda^p$ is a Hermitian Killing form then 
$$d_D^{\star}(\varphi \wedge {J} \varphi)=(d_D^{\star} \varphi) \wedge \mathbb{J} \varphi-\varphi \wedge d_D^{\star} \mathbb{J} \varphi-\frac{2}{p+1} \varphi \b \bpf_D \mathbb{J} \varphi;
$$
\item[(ii)] if $\varphi$ in $\lambda^p$ satisfies 
$$d_D^{\star}(\varphi \wedge \mathbb{J} \varphi)=(d_D^{\star} \varphi) \wedge \mathbb{J} \varphi-\varphi \wedge d_D^{\star} \mathbb{J} \varphi-\frac{2}{p+1} \varphi \b \bpf_D \mathbb{J} \varphi 
$$
and it is everywhere nondegenerate then $\varphi$ is a Hermitian Killing form. 
\end{itemize}
\end{pro}
\begin{proof}
We have 
\begin{equation*}
\begin{split}
e_i \lrcorner D_{e_i}(\varphi \wedge \mathbb{J} \varphi)=&e_i \lrcorner (D_{e_i} \varphi \wedge \mathbb{J} \varphi+\varphi \wedge D_{e_i} \mathbb{J} \varphi)\\
=& (e_i \lrcorner D_{e_i} \varphi) \wedge \mathbb{J} \varphi-D_{e_i} \varphi \wedge (e_i \lrcorner \mathbb{J} \varphi)+(e_i \lrcorner \varphi) \wedge D_{e_i} \mathbb{J} \varphi
-\varphi \wedge (e_i \lrcorner D_{e_i} \mathbb{J} \varphi).
\end{split}
\end{equation*}
After summation, we get 
\begin{equation*}
\begin{split}
-d_D^{\star}(\varphi \wedge \mathbb{J} \varphi)=&-(d_D^{\star} \varphi) \wedge \mathbb{J} \varphi+\varphi \wedge d_D^{\star} \mathbb{J} \varphi\\
&-\sum \limits_{k=1}^{2m} D_{e_i} \varphi \wedge (e_i \lrcorner \mathbb{J} \varphi)+\sum \limits_{k=1}^{2m}(e_i \lrcorner \varphi) \wedge D_{e_i} \mathbb{J} \varphi.
\end{split}
\end{equation*}
But 
\begin{equation*}
\sum \limits_{k=1}^{2m} D_{e_i} \varphi \wedge (e_i \lrcorner \mathbb{J} \varphi)=\sum \limits_{k=1}^{2m} D_{e_i} \varphi \wedge (Je_i \lrcorner \varphi)=-\sum \limits_{k=1}^{2m} D_{Je_i} \varphi \wedge (e_i \lrcorner \varphi)=
-\sum \limits_{k=1}^{2m} (e_i \lrcorner \varphi) \wedge D_{Je_i} \varphi 
\end{equation*}
whence 
%%%%%%%%%%%%%%%%%%%%%%
%\edz{check if can $J$ the frame}
\begin{equation} \label{interequiv}
\begin{split}
-d_D^{\star}(\varphi \wedge \mathbb{J} \varphi)=&-(d_D^{\star} \varphi) \wedge \mathbb{J} \varphi+\varphi \wedge d_D^{\star} \mathbb{J} \varphi\\
& +\sum \limits_{k=1}^{2m} e_i \lrcorner \varphi  \wedge (D_{Je_i} \varphi +D_{e_i}\mathbb{J}\varphi).
\end{split}
\end{equation}
(i) if $\varphi$ is a Hermitian Killing form one uses \eqref{interequiv} and (i) in Proposition \ref{middle} to obtain the conclusion.\\
(ii) in this case, using \eqref{interequiv} again we have that 
$$ \sum \limits_{i=1}^{2m} (e_i \lrcorner \varphi) \wedge \gamma_{e_i}=0
$$
where $\gamma$ in $\lambda^1 \otimes_2 \lambda^p$ is given by 
$$ \gamma_X=D_{JX} \varphi +D_{X}\mathbb{J}\varphi-\frac{2}{p+1} X \lrcorner \bpf_D \mathbb{J} \varphi $$
for all $X$ in $TM$. We define now $\hat{\gamma}$ in $\lambda^{p-1} \otimes \lambda^p$ by 
$$ \hat{\gamma}(X_1, \ldots , X_{p-1})=\gamma_{X_1 \lrcorner \ldots X_{p-1} \varphi}
$$
whenever $X_k, 1 \le k \le p-1$ belong to $TM$. It is easy to verify that $\hat{\gamma}$ is in $\lambda^{p-1} \otimes_1 \lambda^p$, and since $a(\hat{\gamma})=0$ Lemma \ref{l31}, (iii) yields 
$\hat{\gamma}=0$. Because $\varphi$ is everywhere non-degenerate we obtain that $\gamma$ vanishes and we conclude by Proposition \ref{middle}, (i). 
\end{proof}
To end this section let us present another way of characterising Hermitian Killing forms, this time under the form of a product rule with respect to the exterior differential $d_D$.
\begin{pro} \label{schouten}
Let $p \ge 3$ be odd. We have:
\begin{itemize}
\item[(i)] any Hermitian Killing form $\varphi$ in $\lambda^p$ satisfies 
$$ -d_D(\varphi \b\varphi)=(d_D \varphi) \b \varphi+( d_D \mathbb{J} \varphi ) \b \mathbb{J} \varphi-\frac{2}{p+1} J(\varphi \b \bpf_D \mathbb{J}\varphi);
$$
\item[(ii)] if $\varphi$ in $\lambda^p$ satisfies 
$$ -d_D (\varphi \b \varphi)=( d_D \varphi ) \b \varphi+( d_D \mathbb{J} \varphi ) \b \mathbb{J} \varphi-\frac{2}{p+1} J(\varphi \b \bpf_D \mathbb{J}\varphi)
$$
and it is nowhere degenerate, then $\varphi$ is a Hermitian Killing form. 
\end{itemize}
\end{pro} 
\begin{proof} We will mainly use Proposition \ref{id1} and the K\"ahler identities. Indeed, let us observe that 
$$ L^{\star}(\varphi \wedge \mathbb{J} \varphi)=- \varphi \b \varphi.
$$
Using \eqref{kid} it follows that 
$$ d_D L^{\star}(\varphi \wedge \mathbb{J} \varphi)=L^{\star}d_D(\varphi \wedge \mathbb{J} \varphi)+(Jd_D^{\star}J)(\varphi \wedge \mathbb{J} \varphi)
$$
hence 
$$ -d_D(\varphi \b \varphi)=L^{\star}d_D(\varphi \wedge \mathbb{J} \varphi)+Jd_D^{\star}(\varphi \wedge \mathbb{J} \varphi).
$$
It is easy to see that 
\begin{equation*}
\begin{split}
&L^{\star}(d_D \varphi \wedge \mathbb{J} \varphi)=(L^{\star}d_D \varphi) \wedge \mathbb{J} \varphi+(d_D \varphi) \b \varphi\\
&- L^{\star}(\varphi \wedge d_D \mathbb{J}\varphi)=-(L^{\star}d_D \mathbb{J} \varphi) \wedge \varphi+(d_D \mathbb{J} \varphi) \b \mathbb{J} \varphi
\end{split}
\end{equation*}
hence we get further 
\begin{equation*}
\begin{split}
-d_D(\varphi \b \varphi)=&(L^{\star}d_D \varphi) \wedge \mathbb{J} \varphi+(d_D \varphi) \b \varphi\\
&-(L^{\star}d_D \mathbb{J} \varphi) \wedge \varphi+(d_D \mathbb{J} \varphi )\b \mathbb{J} \varphi\\
&+Jd_D^{\star}(\varphi \wedge \mathbb{J} \varphi).
\end{split}
\end{equation*}
Suppose now that $\varphi$ is a Hermitian Killing form. Then, by using (i) in Proposition \ref{id1} we get 
$$ d_D^{\star}(\varphi \wedge \mathbb{J} \varphi)=d_D^{\star} \varphi \wedge \mathbb{J} \varphi-\varphi \wedge d_D^{\star} (\mathbb{J} \varphi)
-\frac{2}{p+1} \varphi \b \bpf_D(\mathbb{J} \varphi)
$$
while the second part of the K\"ahler identities \eqref{kid} provides us with 
\begin{equation*}
\begin{split}
-L^{\star}d_D(\mathbb{J} \varphi)=&d_D^{\star} \varphi\\
L^{\star}d_D \varphi=&d_D^{\star} \mathbb{J} \varphi.
\end{split}
\end{equation*}
The claim in (i) follows now by a simple computation. The converse statement in (ii) is proved by using methods similar to those employed for (ii) in Proposition \ref{id1}, and it is therefore left to the reader.
\end{proof}
%%%%%%%%%%%%%%%%%%%%%%%%% %%%%%%%%%%%%%%%%%%%%%%%%%%%%%%%%%%%%%
\section{The torsion within the $\mathcal{G}_1$ class}
\subsection{$\mathcal{G}_1$-structures} 
In the subsequent $(M^{2m},g,J)$ will be almost-Hermitian, in the class $\mathcal{G}_1$. The connexion $D$ acts on any form $\varphi$ in $\Lambda^{\star}$ according to 
\begin{equation*}
D_X\varphi=\nabla_X\varphi+\frac{1}{2}[X \lrcorner T^D, \varphi]
\end{equation*}
for all $X$ in $TM$. Based on this fact it is straightforward to check that the differential $d_D$ is related to $d$ by 
$$ d_D \varphi=d\varphi-\sum \limits_{k=1}^{2m} (e_k \lrcorner T^D) \wedge (e_k \lrcorner \varphi)
$$
for all $\varphi$ in $\Lambda^{\star}$. This yields, in the particular case of $3$-forms the following comparison formula 
\begin{equation} \label{comp}
d_D\varphi=d \varphi-T^D \bullet \varphi
\end{equation}
for all $\varphi$ in $\Lambda^3$. \\
Let now $R$ be the curvature tensor of the metric $g$ with the convention that 
$R(X,Y)=-\nabla^2_{X,Y}+\nabla^{2}_{Y,X}$ for all $X,Y$ in $TM$. In the same time we consider the curvature tensor $R^D$ of the connexion $D$ where the same 
convention applies. We shall now recall that the first Bianchi identity for the connexion $D$ takes the form:
\begin{lema} \label{b1D} 
$$ b_1(R^D)_X=D_XT^D+\frac{1}{2} X \lrcorner dT^D 
$$
for all $X$ in $TM$.
\end{lema}
\begin{proof} 
It is easy to see, directly from the definition, that the curvature tensors of the Levi-Civita connexion and that of the characteristic connexion are related by 
\begin{equation} \label{curvcomp}
R^D(X,Y)=R(X,Y)-\frac{1}{2} \biggl [ Y \lrcorner D_XT^D- X \lrcorner D_YT^D \biggr ]+\frac{1}{4} \varepsilon^{T^D}(X,Y)
\end{equation}
for all $X,Y$ in $TM$. Here we have set 
$$  \varepsilon^{T^D}(X,Y)=[T^D_X, T^D_Y]-2T^D_{T^D_XY} $$
Since $R$ satisfies the algebraic Bianchi identity and 
$$ b_1(\varepsilon^{T^D})_X=2 X \lrcorner (T^D \b T^D)$$ 
%\footnote{only $[T^D_X, T^D_Y]+2T^D_{T^D_XY} $ is a curvature tensor}
the assertion follows eventually by setting $X=e_k$, taking the exterior product with $e_k$ in \eqref{curvcomp} and summing over $1 \le k \le 2m$. Note that 
in the process one also uses the comparison formula \eqref{comp}.
\end{proof}
Let us gather now some information on the differentials of the components of the torsion form $T^D$.
\begin{pro} \label{diff}
Let $(M^{2m},g,J)$ belong to the class $\mathcal{G}_1$. The following hold:
\begin{itemize}
\item[(i)] $\pf_Dt=0$;
\item[(ii)] $3\pf_D \psi^{+}+2 \bpf_Dt=-4(t \b Jt)-4(Jt \b \psi^{+})_{\lambda^{1,3}}$;
\item[(iii)] $\bpf_D \psi^{+}=-\frac{8}{3}(Jt \b \psi^{+})_{\lambda^4}.$
\end{itemize}
\end{pro}
\begin{proof}
Since $d\omega=2t+3\psi^{+}$ is closed it follows that 
$$ 2dt+3d\psi^{+}=0.
$$
Rewritten by means of the connexion this identity yields, when also using the comparison formula \eqref{comp}
$$ 2(d_Dt+T^D \b t)+3(d_D \psi^{+}+T^D \b \psi^{+})=0.
$$
We now take into account that $T^D=2Jt+\psi^{-}$ to arrive, after expansion of the product and use of Lemma \ref{btype}, (v) at 
$$ 2d_Dt+4Jt \b t+2 \psi^{-} \b t+3d_D \psi^{+}+6 Jt \b \psi^{+}=0.
$$
An elementary computation yields $J(\psi^{-} \b t)=Jt \b \psi^{+}$, in particular 
\begin{equation*}
\begin{split}
&(\psi^{-} \b t)_{\lambda^4}=(Jt \b \psi^{+})_{\lambda^4}\\
&(\psi^{-} \b t)_{\lambda^{1,3}}=-(Jt \b \psi^{+})_{\lambda^{1,3}}.
\end{split}
\end{equation*}
It suffices thus to use Lemma \ref{btype} to obtain the proof of the claims after identifying the various components in the equation above along the bidegree decomposition of $\Lambda^4$. 
\end{proof}
We are now ready to prove the main result of this section.
\begin{teo} \label{main3}
Let $(M^{2m},g,J)$ be almost Hermitian, of type $\mathcal{G}_1$. Then $\psi^{-}$ is a Hermitian Killing form, that is 
$$ D_X \psi^{-}=\biggl [ X \lrcorner (\pf_D \psi^{-}+\frac{1}{4} \bpf_D \psi^{-}) \biggr ]_{\lambda^3}
$$
for all $X$ in $TM$.
\end{teo}
\begin{proof}
Since $D$ is a Hermitian connexion we have $[R^D(X,Y),\omega]=0$, in other words 
$$ \sum \limits_{i=1}^{2m} R^D(X,Y)e_i \wedge Je_i=0
$$
for all $X,Y$ in $TM$. Setting $X=e_k$ and taking the exterior product with $e_k$ we find 
$$ \sum \limits_{i,k=1}^{2m}e_k \wedge R^D(e_k,Y)e_i \wedge Je_i=0.
$$
Obviously, this is equivalent with 
$$ \sum \limits_{i,k=1}^{2m} e_i \lrcorner (e_k \wedge R^D(e_k,Y)) \wedge Je_i= \sum \limits_{i=1}^{2m} R^D(e_i, Y) \wedge Je_i
$$
or further, by using Lemma \ref{b1D}
$$ -\J (D_XT^D+\frac{1}{2} X \lrcorner d_DT^D)=\sum \limits_{i=1}^{2m} R^D(e_i, Y) \wedge Je_i
$$ 
for all $Y$ in $TM$. Since $D$ is Hermitian, we have that $R^D(X,Y)$ is in $\lambda^{1,1}$ hence the right hand side above is in $\lambda^{1,2}$ given that $\Lambda^1 \wedge \lambda^{1,1} \subseteq 
\lambda^{1,2}$. Thus
$$(D_XT^D+\frac{1}{2} X \lrcorner d T^D)_{\lambda^3}=0$$
and it follows that 
\begin{equation} \label{altt}
D_X \psi^{-}=-\frac{1}{2}(X \lrcorner d T^D)_{\lambda^3}
\end{equation}
for all $X$ in $TM$. Therefore, $\psi^{-}$ is a Hermitian Killing form and the claim follows by using Proposition \ref{hk1}. 
\end{proof}
We have the following immediate consequence of the above result.
\begin{teo} \label{par2}
Let $(M^{2m},g,J)$ have totally skew-symmetric Nijenhuis tensor. Then 
\begin{equation*}
DN_J=0 \mbox{\ if and only if \ } dT^D \mbox{\ belongs to \ } \lambda^{2,2}.
\end{equation*} 
In particular $N_J$ is 
parallel if $T^D$ is a closed $3$-form. 
\end{teo}
\begin{proof}
This is a direct consequence of \eqref{altt} and of the properties of the projection on $\lambda^3$. 
\end{proof}
\begin{rema}
It also follows from \eqref{altt} when combined with Theorem \label{main3} that $dT^D$ belongs to $\lambda^{2,2}$ if and only if $d_D\psi^{+}=0$.
\end{rema}
In dimension $6$, the fact that $\psi^{+}$ is a Hermitian Killing form is described by (ii) in Proposition \ref{middle}. As observed in \cite{But2}, it amounts 
then to the local parallelism  of $\psi^{+}$ w.r.t the characteristic connexion, after performing a conformal transformation in order to normalise the 
length of $\psi^{+}$ to a constant. Moreover, J.-B.Butruille shows \cite{But2} that locally 
every $3$-form of suitable algebraic type can be realised as the Nijenhuis form of an almost Hermitian structure of type $\mathcal{G}_1$. In dimension $6$ structures of type $\mathcal{G}_1$ with $DT^D=0$ have been classified 
in \cite{AFS}. Also classification results are avalaible \cite{Al} under the same assumption in the Hermitian case, provided that 
the holonomy of $D$ is contained in $S^1 \times U(m-1)$. 
%%%%%%%%%%%%%%%%%%
\subsection{The class $W_1+W_4$}
In this section we shall record some of the additional features of the geometry  
of  almost Hermitian manifolds in the Gray-Hervella class $W_1+W_4$. It can be described as the subclass of 
$\mathcal{G}_1$ having the property that $t$ has no component on $\Lambda^3_0$, or alternatively 
\begin{equation} \label{NKW}
t=\theta \wedge \omega.
\end{equation}
Here $\theta$ is a $1$-form on $M$, called the Lee form of the almost-Hermitian structure $(g,J)$. It can be recovered directly from 
the K\"ahler form of $(g,J)$ by 
\begin{equation*}
(m-1)J\theta=d^{\star} \omega,
\end{equation*}
in particular we have that $d^{\star}(J\theta)=0$. 
\begin{pro} \label{NKW2}
Let $(M^{2m},g,J), m \ge 3$ be almost-Hermitian in the class $W_1+W_4$. The following hold:
\begin{itemize}
\item[(i)] $\pf_D \theta=0$;
\item[(ii)] $-\pf_D \psi^{+}=2(J\theta \wedge \psi^{-}+\theta \wedge \psi^{+})+\frac{2}{3}(\bpf_D \theta-2 \z \lrcorner \psi^{-}) \wedge \omega$;
\item[(iii)] $ \frac{1}{4} \bpf_D \psi^{+}=\theta \wedge \psi^{+}-J \theta \wedge \psi^{-}$.
\end{itemize}
\end{pro}
\begin{proof}
All claims follows from Proposition \ref{diff}, applied to our present situation, that is $t=\theta \wedge \omega$. Therefore, (i) is immediate from (i) in the 
previously cited Proposition.\\
Through direct computation we find 
\begin{equation*}
\begin{split}
(J\theta \wedge \omega) \b \psi^{+}=&-(\z \lrcorner \psi^{-}) \wedge \omega+3J\theta \wedge \psi^{-} \\
(J\theta \wedge \omega) \b (\theta \wedge \omega)=&0. 
\end{split}
\end{equation*}
Our last two assertions follow now from Proposition \eqref{diff}, after projection of the formulae above onto $\Lambda^4=\lambda^{2,2} \oplus \lambda^{1,3} \oplus \lambda^4$. 
%From \eqref{do} we get, in our setting, 
%$$ d\omega=2\omega \wedge \theta+3\psi^{+}.
%$$ After differentiation, we get 
%\begin{equation*}
%3d\psi^{+}-6\theta \wedge \psi^{+}+2d\theta \wedge \omega=0.
%\end{equation*}
%Using the comparaison formula yields further 
%\begin{equation*}
%3d_D \psi^{+}+2d_D \theta \wedge \omega+3T^D \b \psi^{+}-6 \theta \wedge \psi^{+}+2 (\z \lrcorner T^D) \wedge \omega=0.
%\end{equation*}
%Now, by using \eqref{tD} we have 
%$$ T^D=2J\theta \wedge \omega+\psi^{-}.$$ 
%Through direct computation we find 
%\begin{equation*}
%\begin{split}
%T^D \b \psi^{+}=&2(J\theta \wedge \omega) \b \psi^{+}\\
%=&-2(\z \lrcorner \psi^{-}) \wedge \omega+6 J\theta \wedge \psi^{-} 
%\end{split}
%\end{equation*}
%whence 
%$$ 3d_D \psi^{+}+6(3J\theta \wedge \psi^{-}-\theta \wedge \psi^{+})+2d_D \theta \wedge \omega-4 (\z \lrcorner \psi^{-}) \wedge \omega=0.
%$$
%The claims in (i)-(iii) follow now by projection onto $\Lambda^4=\lambda^{2,2} \oplus \lambda^{1,3} \oplus \lambda^4$. 
\end{proof}
%\begin{coro} Let $(M^{2m},g,J)$ belong to the class $W_1+W_4$. We have:
%\begin{itemize}
%\item[(i)] $d_D^{\star}\psi^{+}=\frac{2(m-2)}{3} \bpf_D \theta-\frac{4(m-5)}{3}\z \lrcorner \psi^{-}$;
%\item[(ii)]
%\end{itemize}
%\end{coro}
%\begin{proof}
%(i) Using \eqref{kid} we find that $L^{\star} \pf_D \psi^{-}=-d_D^{\star} \psi^{+}$, given that $J\psi^{-}=\psi^{+}$ and that $d_D^{\star} \psi^{+}$ belongs to 
%$\lambda^2$. The claim follows now by using (i) in Proposition \ref{NKW2} together with the formula 
%$$ L^{\star}(v \wedge \varphi)=v \wedge L^{\star}(\varphi)+Jv\lrcorner \varphi.
%$$ 
%which holds for any $v$ in $V$ and for any $\varphi$ in $\Lambda^{\star}$. 
%\end{proof}
In dimension $6$ it has been shown in \cite{But2}, \cite{CsIv2} that an almost Hermitian manifold in the class $W_1+W_4$ has closed Lee 
form, that is $d \theta=0$. 
Theorem \ref{main3} combined with Proposition \ref{NKW2} seems a good starting point to investigate up to what extent structures of $W_1+W_4$ are closed in arbitrary dimensions but we shall not pursue in this direction here.
%%%%%%%%%%%%%%%%%%%%%%%%%%%%%%%%%%%%%%%%%5
\subsection{The $\mathfrak{u}(m)$-decomposition of curvature}
In this section we shall investigate, for further use, the splitting of the curvature tensor of the Hermitian connection $D$. Our main goal is to 
identify explicitely the "non-K\"ahler" part of the tensor $R^D$ and to give it an explicit expression in terms of the torsion form $T^D$ of the characteristic connexion. 
Note that for almost quaternionic-Hermitian or $G_2$-structures (not necessarily with skew-symmetric torsion)  similar results have been obtained in \cite{CaSwcurv, CsIv}.
We will mainly use that 
\begin{equation} \label{skewpart}
\begin{split}
R^D(X,Y,Z,U)-R^D(Z,U,X,Y)=& -\frac{1}{2}\biggl [ (D_XT^D)_Y(Z,U)-(D_YT^D)_X(Z,U) \biggr ]\\
+&\frac{1}{2}\biggl [  (D_ZT^D)_U(X,Y)-(D_ZT^D)_U(X,Y) \biggr ]
\end{split}
\end{equation}
holds for all $X,Y,Z,U$ in $TM$, as it easily follows from the difference formula \eqref{curvcomp}, after verifying that $\varepsilon^{T^D}$ belongs to $S^2(\Lambda^2)$. 
\begin{teo} \label{curvfin1}
Let $(M^{2m},g,J)$ be almost Hermitian of type $\mathcal{G}_1$. We have:
\begin{equation*}
R^D=R^K+\hat{\Omega}+\frac{1}{2}R^a+R^m
\end{equation*}
where $R^K$ belongs to $\mathcal{K}(\mathfrak{u}(m))$. Moreover:
\begin{itemize}
\item[(i)] $R^a$ belongs to $\Lambda^{2}(\lambda^{1,1})$ and satisfies the Bianchi identity  
\begin{equation*}
\frac{1}{2}(b_1 R^a)_X=D_X(Jt)-D_{JX}t-\frac{1}{2}X \lrcorner \pf_D(Jt)\\
\end{equation*}
for all $X$ in $TM$.
\item[(ii)] the Bianchi identity for $R^m$ in $\lambda^2 \otimes \lambda^{1,1}$ is 
\begin{equation*}
(b_1 R^m)_X=D_X(Jt)+D_{JX}t+\frac{1}{4}(JX \lrcorner \J \pf_D \psi^{-}-2 X \lrcorner \pf_D \psi^{-})\\\
\end{equation*}
for all $X$ in $TM$.
\item[(iii)] $\Omega$ in $\lambda^{2,2}$ is given by 
\begin{equation*}
\Omega=\frac{3}{2} \pf_D(Jt)+2(Jt \b Jt)_{\lambda^{2,2}}+\frac{1}{2} \psi^{-} \b \psi^{-}
\end{equation*}
\end{itemize}
\end{teo}
\begin{proof}
This is a direct application of Theorem \ref{curvfin0} which can be used, as it follows from \eqref{skewpart}, for the tensor 
$\RR=R^D$ and 
$$ \gamma=-\frac{1}{2}DT^D=-D(Jt)-\frac{1}{2}D \psi^{-}.
$$ 
Then we have $\gamma^{1,2}=-D(Jt)$ and $\gamma^3=-\frac{1}{2}D\psi^{-}$, and in order to prove our claims we only need to determine the various antisymmetrisations 
of $\gamma^{1,2}$. Also note that in this context
$$ T=-\frac{1}{2}d_DT^D=-d_D(Jt)-\frac{1}{2}d_D \psi^{-},
$$
in particular $T^{2,2}=-\pf_D (Jt)$.\\
(i) Now, it is easy to see that 
\begin{equation*}
\begin{split}
a(\gamma^{1,2})=&-d_D(Jt)\\
a(J \gamma^{1,2})=& d_Dt \\
a^c(\gamma^{1,2})=&-Jd_Dt\\
a^c(J\gamma^{1,2})=&-Jd_DJt.
\end{split}
\end{equation*}
Hence,
\begin{equation*}
\begin{split}
A_1=&-d_D(Jt)-Jd_D(Jt)-4T^{2,2}=-d_D(Jt)-Jd_D(Jt)+4\pf_D(Jt)\\
=&-2\pf_D(Jt)+4\pf_D(Jt)=2 \pf_D(Jt)
\end{split}
\end{equation*}
and 
$$A_2=-Jd_Dt-d_Dt=-2\pf_Dt=0
$$
by (i) of Proposition \ref{diff}.\\
(ii) Since $T^{1,3}=-\bpf_D(Jt)-\frac{1}{2} \pf_D \psi^{-}$ we have 
\begin{equation*}
\begin{split}
A_3=&-d_D(Jt)+Jd_DJt+2\bpf_D(Jt)+\pf_D \psi^{-}\\
=&\pf_D \psi^{-}
\end{split}
\end{equation*}
and 
\begin{equation*}
\begin{split}
A_4=& d_Dt-Jd_Dt+\J\bpf_D(Jt)+\frac{1}{2} \J \pf_D \psi^{-}\\
=&\frac{1}{2} \J \pf_D \psi^{-}
\end{split}
\end{equation*}
after making use of \eqref{dtype}. The claim in (iii) now follows.\\
(iii) We apply the Bianchi operator to 
$$ R^D=R^K+\hat{\Omega}+R^a+R^m $$
and find after making use of (ii), (iii) and of \eqref{b1f} that 
\begin{equation*}
b_1(R^D)_X=X \lrcorner \Omega+2D_X(Jt)-\frac{1}{2} X \lrcorner \pf_D(Jt)+\frac{1}{4}(JX \lrcorner \J \pf_D \psi^{-}-2 X \lrcorner \pf_D \psi^{-})
\end{equation*}
for all $X$ in $TM$. Plugging into this the Bianchi identity for $D$ in Lemma \ref{b1D}, combined with the fact that $\psi^{-}$ is a Hermitian Killing form as asserted in Theorem \ref{main3} yields further 
\begin{equation*}
\begin{split}
\frac{1}{2}(X \lrcorner dT^D)_{\lambda^{1,2}}=&X \lrcorner \biggl [ \Omega-\frac{1}{2} \pf_D(Jt) \biggr ]+\frac{1}{4}(JX \lrcorner \J \pf_D \psi^{-}-2 X \lrcorner \pf_D \psi^{-})\\
=&X \lrcorner \biggl [ \Omega-\frac{1}{2} \pf_D(Jt) \biggr ]-(X \lrcorner \pf_D \psi^{+})_{\lambda^{1,2}}
\end{split}
\end{equation*}
for all $X$ in $TM$. After identifying components on $\lambda^1 \otimes_{1} \lambda^{1,2}$ and $\lambda^1 \otimes_2 \lambda^{1,2}$ respectively (see also Remark 2.1 for the definition of these spaces) we find that 
$$ \frac{1}{2}(dT^D)_{\lambda^{2,2}}= \Omega-\frac{1}{2} \pf_D(Jt).
$$
Now $dT^D=d_DT^D+T^D \b T^D$ by the comparaison formula \eqref{comp} thus by projection on $\lambda^{2,2}$ we find 
\begin{equation*}
\begin{split}
(dT^D)_{\lambda^{2,2}}=&2\pf_D(Jt)+(T^D \bullet T^D)_{\lambda^{2,2}}\\
=&2\pf_D(Jt)+4(Jt \b Jt)_{\lambda^{2,2}}+\psi^{-} \b \psi^{-}
\end{split}
\end{equation*}
after expansion of the product and use of Lemma \ref{btype}, (iii). The claim in (iii) follows now immediately. 
\end{proof}
\begin{rema}
The curvature decomposition in Theorem \ref{curvfin1} can be still refined, given that the $U(m)$-modules 
$$\K(\mathfrak{u}(m)), \lambda^{2,2}, \Lambda^2(\lambda^{1,1}) \ \mbox{and} \ \lambda^2 \otimes \lambda^{1,1}$$ 
are not irreducible. Altough this is an algebraically simple procedure, the computations at the level of the derivative $DT^D$ of the torsion form become 
somewhat involved and will not be presented here. We just illustrate the situation in the simpler case of $W_1+W_4$ below.
\end{rema}
Theorem \ref{curvfin1} has various applications, a class of which consists in giving 
torsion interpretation of curvature conditions imposed on the tensors $R^D$ or $R$. As an example in this direction we have the following:
\begin{pro} \label{herm}
Let $(M^{2m},g,J)$ belong to the class $\mathcal{G}_1$. The curvature tensor $R^D$ is Hermitian, that is 
$$ R^D(JX, JY)=R^D(X,Y)
$$
for all $X,Y$ in $TM$, if and only if 
\begin{equation*}
D_{JX}t+D_X(Jt)=\frac{2}{3} (X \lrcorner \bpf_D(Jt))_{\lambda^{1,2}}
\end{equation*}
for all $X$ in $TM$ and 
$$ 2\bpf_D(Jt)=3\pf_D\psi^{-}.
$$
\end{pro}
\begin{proof}
By Theorem \ref{curvfin1} the tensor $R^D$ is Hermitian if and only if $R^m=0$ which in turn happens if and only if $b_1R^m=0$. Therefore, by (ii) in Theorem \ref{curvfin1} the curvature of $D$ is Hermitian 
if and only if 
$$D_X(Jt)+D_{JX}t+\frac{1}{4}(JX \lrcorner \J \pf_D \psi^{-}-2 X \lrcorner \pf_D \psi^{-})=0$$
for all $X$ in $TM$. Taking the alternation above we find
$$ 2\bpf_D(Jt)-3\pf_D\psi^{-}=0.
$$
We conclude by recalling that 
\begin{equation} \label{proj12}
(X \lrcorner A)_{\lambda^{1,2}}=\frac{1}{4}(2X \lrcorner A-JX \lrcorner \J A)
\end{equation}
for all $X$ in $TM$, where $A$ belongs to $\lambda^{1,3}$.
\end{proof}
\begin{rema} \label{rcurvhermalg}
In the case when the curvature tensor of the connection $D$ is Hermitian the fact that $2\bpf_D(Jt)=3\pf_D\psi^{-}$ combined with (ii) in Proposition \ref{diff} yields further
$$ \bpf_D(Jt)=-(t \b Jt)-(Jt \b \psi^{+})_{\lambda^{1,3}}.
$$
\end{rema}
For the subclass $W_1 \oplus W_4 \subseteq \mathcal{G}_1$ more information on the curvature tensor of the connection $D$ is available and it will actually turn out that the components 
$R^a$ and $R^m$ have simple algebraic expressions. We define 
$$ S^{2,-}(M)=\{S \in S^2(M) : SJ+JS=0\}.
$$
This is embedded in $\lambda^2 \otimes \lambda^{1,1}$ via $S \mapsto \mathring{S}$ where 
$$ \mathring{S}(X,Y)=\frac{1}{2}(SJX\wedge Y+X \wedge SJY+SX \wedge JY+JX \wedge SY)
$$
for all $X,Y$ in $TM$. One verifies that 
\begin{equation} \label{ering}
(b_1 \mathring{S})_X=SX \wedge \omega
\end{equation}
for all $X$ in $TM$. We also have an embedding $\lambda^{1,3} \hookrightarrow \lambda^2 \otimes \lambda^{1,1}$ given by $\Omega \mapsto \tilde{\Omega}$ where we define 
\begin{equation*}
\tilde{\Omega}(X,Y)=\frac{1}{4}(\Omega(JX,JY)-\Omega(X,Y))
\end{equation*}
for all $X,Y$ in $TM$. Elementary considerations ensure that this is well defined and subject to 
\begin{equation} \label{etilde}
(b_1 \tilde{\Omega})_X=(X \lrcorner \Omega)_{\lambda^{1,2}}
\end{equation}
whenever $X$ belongs to $TM$. As a last piece of notation let the symmetrised action of $D$ on $1$-forms be defined by 
$$(\mathring{D}\alpha)(X,Y)=\frac{1}{2}((D_X\alpha)Y+(D_Y\alpha)X)$$ 
whenever $\alpha$ is in $\lambda^1$ and $X,Y$ belong to $TM$.
\begin{pro} \label{CNKW} Let $(M^{2m},g,J)$ belong to the class $W_1+W_4$. If $\theta$ denotes the Lee form of $(g,J)$ the folowing hold:
\begin{itemize}
\item[(i)] $\Omega=\biggl [ \frac{3}{2} \pf_D(J\theta)+2\vert \theta \vert^2 \omega-4 \theta \wedge J\theta  \biggr ] \wedge \omega+\frac{1}{2} \psi^{-} \b \psi^{-}$;
\item[(ii)] 
$$R^a=\pf_D(J\theta) \otimes \omega-\omega \otimes \pf_D(J\theta);$$
\item[(iii)] 
\begin{equation*}
R^m=-\bpf_D(J\theta) \otimes \omega+\mathring{S}_{\theta}-\widetilde{\pf_D \psi^{-}}
\end{equation*}
where $S_{\theta}$ in $S^{2,-}(M)$ is defined by $S_{\theta}=(1-J)\mathring{D}(J\theta)$. 
\end{itemize}
\end{pro}
\begin{proof} Recall that for class $W_1+W_4$ we have that $t=\theta \wedge \omega$. Taking this into account we will apply now Theorem \ref{curvfin1}.\\
(i) follows from (iii) inTheorem \ref{curvfin1} when observing that  
$$ (Jt \b Jt)_{\lambda^{2,2}}=\vert \theta \vert^2 \omega \wedge \omega-2\omega \wedge \theta \wedge J\theta.
$$
(ii) First of all we note that 
$$ D_X(J\theta)-D_{JX}\theta=X \lrcorner \pf_D(J\theta)$$
for all $X$ in $TM$. The quickest way to see this is to observe that the tensor $q$ in $\lambda^1 \otimes \lambda^1$ defined by 
$q(X)= D_X(J\theta)-D_{JX}\theta-X \lrcorner \pf_D(J\theta)$
for all $X$ in $TM$ is actually in $\lambda^1 \otimes_1 \lambda^1$ and satisfies $a(q)=a^c(q)=0$, by using also that 
$\pf_D\theta=0$ (cf. Proposition \ref{NKW2}, (i)). Our claim now follows by observing that (iii) in Lemma \ref{l31} continues to hold when $p=q=1$.\\
(ii) Using now (i) in Theorem \ref{curvfin1} we get 
$$ (b_1R^a)_X=(X \lrcorner \pf_D(J\theta) \wedge \omega-\pf_D(J\theta) \wedge (X \lrcorner \omega)
$$
for all $X$ in $TM$.\\
Now for any $\varphi$ in $\lambda^{1,1}$ we consider the element $\omega \otimes \varphi-\varphi \otimes \omega$ in $\Lambda^{2}(\lambda^{1,1})$, which is explicitely given by 
$$ (\omega \otimes \varphi-\varphi \otimes \omega)(X,Y)=\omega(X,Y) \varphi-\varphi(X,Y) \omega$$
for all $X,Y$ in $TM$. A straightforward computation 
following the definitions yields
$$ b_1(\pf_D(J\theta) \otimes \omega-\omega \otimes \pf_D(J\theta))_X=(X \lrcorner \pf_D(J\theta)) \wedge \omega-\pf_D(J\theta) \wedge (X \lrcorner \omega)
$$
for all $X$ in $TM$, and we conclude by using the injectivity of the Bianchi map $b_1 : \Lambda^2(\Lambda^2) \to \Lambda^1 \otimes \Lambda^3$.\\
(iii) By (iii) in Theorem \ref{curvfin1}  we have 
\begin{equation*}
(b_1R^m)_X=\biggl [ D_X(J\theta)+D_{JX} \theta \biggr ] \wedge \omega-(X \lrcorner \pf_D \psi^{-})_{\lambda^{1,2}}
\end{equation*}
after also making use of \eqref{proj12}. Since 
$$ D_X(J\theta)+D_{JX} \theta=X \lrcorner \bpf_D(J\theta)+S_{\theta}X$$
it follows by means of \eqref{ering}, \eqref{etilde}, that $R ^m+\bpf_D(J\theta) \otimes \omega-\mathring{S}_{\theta}+\widetilde{\pf_D \psi^{-}}
$ belongs to $Ker(b_1) \cap (\lambda^2 \otimes \lambda^{1,1})$. It therefore vanishes and the proof is finished.
\end{proof}
Finally, we remark that $R^m$ can be given a more detailed expression by taking into account (ii) of Proposition \ref{NKW2}. 

%%%%%%%%%%%%%%%%%%%%%%%%%%%%%%%%%%%%%%%%%%%
\section{Nearly-K\"ahler geometry}
We shall begin here our survey of nearly-K\"ahler geometry. This is the class of almost Hermitian structures introduced by the following.
\begin{defi} \label{nk1}
Let $(M^{2m},g,J)$ be an almost Hermitian manifold. It is called nearly-K\"ahler, shortly (NK) if and only if 
$$ (\nabla_XJ)X=0
$$
whenever $X$ belongs to $TM$. 
\end{defi}
This can be easily rephrased to say that an almost Hermitian structure $(g,J)$ is nearly-K\"ahler if and only its K\"ahler form is subject to 
\begin{equation} \label{kilf}
\nabla_X\omega=\frac{1}{3}X \lrcorner d \omega
\end{equation}
for all $X$ in $TM$. In other words, the K\"ahler form of any NK-structure is 
a Killing form and conversely U.Semmelmann \cite{uwe1} shows that any almost Hermitian structure with this property must be nearly-K\"ahler. 
Therefore, nearly-K\"ahler structures belong to the class $W_1$ in the sense that the $3$-form $t$ in Proposition \ref{st1} vanishes identically. \\
From now on we shall work on a given nearly-K\"ahler manifold $(M^{2m},g,J)$. For many of the properties of an NK-structure are best expressed by means of its first canonical Hermitian connection, 
it is important to note then the coincidence of the canonical connection and $D$, that is 
$$ D=\bnabla.
$$
For this reason, the torsion tensor $T^D$ will be denoted simply by $T$ in what follows. It is given by 
$$ T^D=\psi^{-}
$$
and hence belongs to $\lambda^3$. Therefore, in dimensions $4$ NK-structures are K\"ahler and we shall assume from now on that $m \ge 3$. 
The Nijenhuis tensor of the almost complex $J$ is computed from \eqref{nlc} by 
$$ N^J=-4\psi^{-}.
$$ 
As an immediate consequence one infers that the almost complex structure of an NK-manifold is integrable if and only if 
the structure is actually a K\"ahler one. This observation motivates the following
\begin{defi} \label{strict} 
Let $(M^{2m},g,J)$ be a NK-structure. It is called strict iff $\nabla_XJ=0$ implies that $X=0$ for all $X$ in $TM$.
\end{defi}
Therefore the non-degeneracy of any of the forms $\psi^{+}$ and $\psi^{-}$ is equivalent with the strictness of the NK-structure 
$(g,J)$. 
An important property of NK-structures is contained in the following result. 
\begin{teo} \label{NK1}
Let $(M^{2m},g,J)$ be a nearly-K\"ahler manifold. Then 
$$ \bnabla \psi^{\pm}=0
$$
in other words the Nijenhuis tensor of $J$ is parallel w.r.t the canonical Hermitian connection. 
\end{teo}
This has been proved first by Kirichenko \cite{kiri} and a short proof can be found in \cite{BeMor}. It also follows from 
our Theorem \ref{main3}, which is unifying this type of property in the class $\mathcal{G}_1$. \\
Let us introduce now the symmetric tensor $r$ in $S^2M$ given by 
\begin{equation*}
\langle rX,Y \rangle=\langle X \lrcorner \psi^{+}, Y \lrcorner \psi^{+}\rangle 
\end{equation*}
for all $X,Y$ in $TM$. It is easily seen to be $J$-invariant and if $(M^{2m},g,J)$ is strict then $r$ is non-degenerate. Moreover, from Theorem \ref{NK1} it also 
follows that 
$$ \bnabla r=0.
$$
\begin{coro} \cite{Nagy1} Any NK-manifold is locally the product of a K\"ahler manifold and a strict nearly-K\"ahler one.
\end{coro}
\begin{proof}
Consider the $\bnabla$-parallel distribution $\V:=\{V \in TM : \psi^{+}_V=0\}$. Since the torsion of $\bnabla$ vanishes in direction of $\V$ the latter must be parallel for the Levi-Civita 
and the result follows by using the deRham splitting theorem.
\end{proof}
It follows that locally and also globally if our original manifold is simply connected we can restrict to the study of strict nearly-K\"ahler structures, for short (SNK).
Note however that in dimension $6$, any NK-structure which is not K\"ahler must be strict. \\
The parallelism of the torsion tensor w.r.t $\bnabla$ is eventually reflected in the properties its curvature tensor. Indeed
\begin{pro} \cite{g1,g3} \label{curv}
The following hold:
\begin{itemize}
\item[(i)] $\overline{R}(X,Y,Z,U)=\overline{R}(Z,U,X,Y)$ for all $X,Y,Z,U$ in $TM$;
\item[(ii)] $\overline{R}(JX,JY)=\overline{R}(X,Y)$ for all $X,Y$ in $TM$;
\item[(iii)] 
\begin{equation*}
\overline{R}(X,Y)Z+\overline{R}(Y,X)Z+\overline{R}(Z,X)Y=[\psi^{-}_X, \psi^{-}_Y]Z-\psi^{-}_{\psi^{-}_XY}Z
\end{equation*}
for all $X,Y,Z$ in $TM$. 
\end{itemize}
\end{pro}
\begin{proof} (i) is immediate from \eqref{skewpart} and the parallelism of the torsion, whilst (ii) follows from (i) and the fact that $\overline{R}$ belongs to $\Lambda^2 \otimes \lambda^{1,1}$. The last claim follows 
for instance from Lemma \ref{b1D}
\end{proof}
\begin{rema} Proposition \ref{curv} continues to hold when the nearly-K\"ahler metric is allowed to have signature. This has been exploited in \cite{coscha} to 
classify NK-structures compatible with a flat pseudo-Riemannian metric.
\end{rema}
To establish first order properties of SNK-structures we will have a look at the Ricci tensor of such metrics. The Hermitian Ricci tensor of $(g,J)$ is defined by 
\begin{equation*}
\langle \overline{Ric}X,Y\rangle=\sum \limits_{i=1}^{2m}\overline{R}(X,e_i,Y,e_i)
\end{equation*}
for all $X,Y$ in $TM$, where $\{e_i, 1 \le i \le 2m\}$ is some local orthonormal frame. by making use of Proposition \ref{curv} we find that $\overline{Ric}$ is actually symmetric and 
$J$-invariant. The Hermitian Ricci tensor is related to the usual Riemannian one by 
\begin{equation} \label{Ricibric}
\overline{Ric}=Ric-\frac{3}{4}r
\end{equation}
as implied by the general curvature comparison formula \eqref{curvcomp}.
\begin{teo} \cite{Nagy1}
Let $(M^{2m},g,J)$ be a strict nearly-K\"ahler manifold. The following hold:
\begin{itemize}
\item[(i)] the Ricci tensor of the metric $g$ is parallel w.r.t  $\bnabla$, that is $\bnabla Ric=0$;
\item[(ii)] $Ric$ is positive definite.
\end{itemize}
\end{teo}
\begin{proof}
The proof of both relies on the explicit computation of the Ricci tensor of the metric $g$. From the parallelism of $\psi^{+}$, after derivation and use of the Ricci identity for the connection with torsion $\bnabla$ 
we find 
\begin{equation} \label{isoobs}
[\overline{R}(X,Y), \psi^{+}]=0
\end{equation}
for all $X,Y$ in $TM$. In a local orthonormal frame $\{ e_i, 1 \le i \le 2m\}$ this reads 
$$ \sum \limits_{i=1}^{2m} \overline{R}(X,Y)e_i \wedge \psi^{+}_{e_i}=0
$$
for all $X,Y$ in $TM$. We now set $Y=e_k$, take the interior product with $e_k$ above to find, after summation over $1 \le k \le 2m$ and some straightforward manipulations that 
\begin{equation*}
\psi^{+}_{\overline{Ric}X}=\sum \limits_{1 \le k,i \le 2m} \overline{R}(X,e_k)e_i \wedge \psi^{+}_{e_i, e_k}.
\end{equation*}
Since $\psi^{+}$ is a form, the sum in the right hand side equals 
\begin{equation*}
\begin{split}
&\frac{1}{2}\sum \limits_{1 \le k,i \le 2m}( \overline{R}(X,e_k)e_i-\overline{R}(X,e_i)e_k ) \wedge \psi^{+}_{e_i, e_k}\\
=&\frac{1}{2} \sum  \limits_{1 \le k,i \le 2m} \biggl ( [\psi^{+}_X, \psi^{+}_{e_k}]e_i-\psi^{+}_{\psi^{+}_Xe_k }e_i \biggr ) \wedge \psi^{+}_{e_i, e_k}\\
+&\frac{1}{2}\sum \limits_{1 \le k,i \le 2m}( \overline{R}(e_i,e_k)X \wedge  \psi^{+}_{e_i, e_k}
\end{split}
\end{equation*}
after using the Bianchi identity for $\overline{R}$. Now the last sum vanishes since $\psi^{+}$ is $J$-anti-invariant in the first arguments whereas $\overline{R}$ is $J$-invariant 
and our frame can be chosen to be Hermitian. After neglecting terms which are $J$-invariant in $e_i,e_k$ in the algebraic sum above we end up with 
\begin{equation*}
\psi^{+}_{\overline{Ric}X}=\frac{1}{2}\sum  \limits_{1 \le k,i \le 2m} \psi^{+}_X \psi^{+}_{e_k}e_i \wedge \psi^{+}_{e_i, e_k}
\end{equation*}
for all $X$ in $TM$. Using now the definition of the tensor $r$ a straightforward manipulation yields 
\begin{equation} \label{ricf}
\psi^{+}_{\overline{Ric}X}=-\frac{1}{2} \sum \limits_{k=1}^{2m} \psi^{+}_Xe_k \wedge re_k
\end{equation}
for all $X$ in $TM$. By derivation and using the parallelism of $\psi^{+}$ it follows that 
$$ \psi^{+}_{(\bnabla_X \overline{Ric})Y}=0
$$
for all $X,Y$ in $TM$. (i) follows now from the fact that $\psi^{+}$ is nondegenerate and the comparison fact in \eqref{Ricibric}. For the claim in (ii) we refer the reader to \cite{Nagy1}. 
\end{proof}
Using Myer's theorem it follows from the above that complete SNK-manifolds must be compact with finite fundamental group. Therefore, in the compact 
case one can restrict attention, up to a finite cover, to simply connected nearly-K\"ahler manifolds.\par
Another important object is the first Chern form of the almost Hermitian structure $(g,J)$ defined by 
$$ 8\pi \gamma_1(X,Y)=\sum \limits_{i=1}^{2m}\overline{R}(X,Y,e_i, J_ei)
$$
for all $X,Y$ in $TM$. Since $\overline{R}$ is $J$-invariant we have that $\gamma_1$ belongs to $\lambda^{1,1}$ so one can write $4\pi \gamma_1=\langle CJ \cdot, \cdot \rangle$ for some 
$C$ in $S^2(TM)$ such that $CJ=JC$. A straightforward computation using the first Bianchi identity yields the relation 
\begin{equation*}
C=\overline{Ric}-r
\end{equation*}
hence $C$ must be parallel w.r.t. the canonical connexion, that is 
$$ \bnabla C=0
$$
by means of Theorem 4.2. Since $\gamma_1$ is a closed form, as it follows from the second Bianchi identity for $\bnabla$, this results in the algebraic obstruction 
\begin{equation} \label{chern}
-C(\psi^{+}_XY)=\psi^{+}_XCY+\psi^{+}_{CX}Y
\end{equation}
for all $X,Y$ in $TM$. Note that this can be given a direct algebraic proof by observing that $[\gamma_1,\psi^{+}]=0$, as it follows from \ref{isoobs} when taking into account that $\overline{R}$ belongs to 
$S^2(\Lambda^2)$. 
\subsection{The irreducible case}
To obtain classification results for strict nearly-K\"ahler manifolds we shall start from examining the holonomy representation of the canonical 
Hermitian connection. At a point $x$ of $M$ where $(M^{2m},g,J)$ is some SNK-manifold this is the representation 
$$ Hol_x(\bnabla) : T_xM \to T_xM
$$
obtained by parallel transport w.r.t $\bnabla$ along loops about $x$. The holonomy representation is Hermitian, for $\bnabla$ is a Hermitian connection 
and this gives rise to two different notions of irreducibility as indicated by the following well known result from representation theory.
\begin{pro} \label{irep}
Let $(V^{2m},g,J)$ be a Hermitian vector space and let $(G,V)$ be a Hermitian representation of some group $G$. If we write 
$V^{\mathbb{C}}$ for the complex vector space obtain from $V$ by setting $iv=Jv$ for all $v$ in $V$ the following cases can occur:
\begin{itemize}
\item[(i)] $(G,V)$ is irreducible;
\item[(ii)] $(G,V^{\mathbb{C}})$ is irreducible but not $(G,V)$. In this case $V$ splits orthogonally as $V=L \oplus JL$ for some 
$G$-invariant subspace $L$ of $V$;
\item[(iii)] $(G,V^{\mathbb{C}})$ is reducible. 
\end{itemize}
\end{pro}
In this section we shall deal with the instances when the holonomy representation of $\bnabla$ is irreducible in the sense of (i) or (ii) in the
Proposition \ref{irep}. The first is actually covered by the following powerful result of R.Cleyton and A.Swann. 
\begin{teo} \cite{cs} \label{cswann}
Let $(N^{n},g)$ be Riemannian such that there exists a metric connection $D$ such that the following hold:
\begin{itemize}
\item[(i)] the torsion tensor $T$ of $D$ belongs to $\Lambda^3$;
\item[(ii)] $DT=0$ and $T$ does not vanish identically.
\end{itemize}
If the holonomy representation of $D$ is irreducible then $D$ is an 
Ambrose-Singer connection in the sense that $DR^{D}=0$ where $R^D$ denotes the curvature tensor of the connection $D$, provided that $n \neq 6,7$. 
\end{teo}
The two exceptions in the result above correspond actually to nearly-parallel $G_2$-structures in dimensions $7$ (see \cite{FKMS} for an account) and SNK-structures in dimension $6$.\\
In the situation in the Theorem above $(N^n,g)$ is a locally homogeneous space (see \cite{TriVanhecke2}) for more details). Theorem \ref{cswann} is proved 
by making use of general structure results on Berger algebras and formal curvature tensors spaces, for irreducible representations of compact Lie 
algebras (see also \cite{MerkSwach} for the non-compact case). We can now make the following.
\begin{teo} \label{irednk} Let $(M^{2m},g,J)$ be a strict nearly-K\"ahler manifold. Then either:
\begin{itemize}
\item[(i)] $\bnabla$ is an Ambrose-Singer connection;
\item[(ii)] $m=3$;
\item[] or 
\item[(ii)] the holonomy representation of $\bnabla$ is reducible over $\mathbb{C}$.
\end{itemize}
\end{teo}
\begin{proof}
By making use of the Theorem \ref{cswann} we see that the case (i) in Proposition \ref{irep} corresponds to (i) in our statement. To finish 
the proof let us suppose that at some point $x$ of $M$ we have an orthogonal splitting 
$$ T_xM=L_x \oplus JL_x
$$
for some $Hol_x(\bnabla)$-invariant subspace of $T_xM$. Using parallel transport $L_x$ extends to a $\bnabla$-parallel distribution $L$ of 
$TM$ such that 
$$ TM=L \oplus JL.
$$
It follows that $\overline{R}(L,JL)=0$ by also using Proposition \ref{curv}, (ii). Since $\overline{R}$ belongs to $S^2(\lambda^{1,1})$, this means that 
$\overline{R}$ is an algebraic expression in $\psi^{+}$, which can be explicitly determined from the first Bianchi identity (see \cite{Nagy2} for details).  It follows that  
that $\bnabla \ \overline{R}=0$ and then $\bnabla$ is an Ambrose-Singer connection. 
\end{proof}
\section{When the holonomy is reducible}
In this section we present classification results for SNK-structures in the case when the holonomy representation of $\bnabla$ is complex reducible. We begin by setting up some 
terminology which is aimed to gain some understanding concerning the relation between the algebraic properties of the torsion form of a SNK-structure and the geometry 
of the underlying Riemmanian manifold.
\subsection{Nearly-K\"ahler holonomy systems}
We start by the following definition which extracts the more peculiar facts from NK-geometry which relate to the holonomy of the canonical Hermitian connection.
\begin{defi} \label{nkhsys}
A nearly-K\"ahler holonomy system $(V^{2m},g,J, \psi^{+},R)$ is constituted of the following data:
\begin{itemize}
\item[(i)] a Hermitian vector space $(V^{2m},g,J)$;
\item[(ii)] a non-zero $3$-form $\psi^{+}$ in $\lambda^3$;
\item[(iii)] a tensor $R$ in $\lambda^{1,1} \otimes \lambda^{1,1}$ of the form 
\begin{equation*}
R=R^K+\hat{\Omega}
\end{equation*}
where $R^K$ in $\mathcal{K}(\mathfrak{u}(m))$ is such that 
\begin{equation} \label{mm1}
[R(x,y), \psi^{+}]=0
\end{equation}
holds for any $x,y$ in $V$. Moreover the form $\Omega$ in $\lambda^{2,2}$ must be given by $\Omega=\frac{1}{2}\psi^{+} \b \psi^{+}$. 
\end{itemize}
\end{defi}
A given NK-holonomy system $(V^{2m},g,J, \psi^{+}, R)$ is called strict if and only if the form $\psi^{+}$ is non-degenerate. In what follows we shall work 
on a given SNK-holonomy system to be denoted by $(V^{2m},g,J, \psi^{+}, R)$. As with Riemannian holonomy systems, the liaison with the holonomy algebra of the canonical connection of a geometric 
SNK-structure is through the Lie subalgebra 
$$ \mathfrak{h}:=Lie \{R(x,y) : x,y \in V \}
$$
of $\lambda^{1,1}$. 
\begin{defi} \label{cred}
An SNK-holonomy system $(V^{2m},g,J, \psi^{+}, R)$ is complex reducible if the representation $(\mathfrak{h},V^{\mathbb{C}})$ is reducible.  
\end{defi}
Another object of relevance here is the isotropy algebra $\g$ of the form $\psi^{+}$ defined by 
\begin{equation*}
\g=\{ \alpha \in \Lambda^2 : [\alpha, \psi^{+}]=0 \}.
\end{equation*}
It is easy to see, starting from \eqref{comrule} and then using an invariance argument together with the non-degeneracy of $\psi^{+}$ that 
$$\g \subseteq \lambda^{1,1}.$$
Moreover, the condition \eqref{mm1} in Definition \ref{nkhsys} reads $\mathfrak{h} \subseteq \g$. \\
We are interested here 
in the structure of complex invariant subspaces of the metric representation $(\mathfrak{h},V)$. The following definition singles out three main classes of subspaces of relevance for our situation.
\begin{defi} \label{subspc}
A proper, $J$-invariant subspace $\V$ of $V$ is said to be (w.r.t. $\psi^{+}$):
\begin{itemize}
\item[(i)] isotropic if $\psi^{+}(\V,\V) \subseteq \V$;
\item[(ii)] null if $\psi^{+}(\V,\V)=0$;
\item[(iii)] special if it is null and $\psi^{+}(H,H)=\V$, where $H$ is the orthogonal complement of $\V$ in $V$. 
\end{itemize}
\end{defi}
\begin{rema} Any $2$-dimensional, $J$-invariant subspace of $V$ is null w.r.t. $\psi^{+}$. Moreover, in dimension $6$, any two dimensional $J$-invariant subspace is 
special w.r.t. $\psi^{+}$. However, we are interested here in isotropic or special subspaces which are invariant w.r.t a particular Lie group or Lie algebra. 
\end{rema}
A useful criterion to prove that a subspace is special is the following. 
\begin{lema} \label{special}
Let $(V^{2m},g, J)$ be a Hermitian vector space and let $\psi^{+}$ be non-degenerate in $\lambda^3$.  If $\V  \subseteq V$ is null w.r.t. $\psi^{+}$ and such that 
$\psi^{+}(H,H) \subseteq \V $, then:
\begin{itemize}
\item[(i)] $\V$ is special w.r.t. $\psi^{+}$;
\item[(ii)] $\psi^{+}(\V,H)=H$.
\end{itemize}
\end{lema}
\begin{proof}
(i) First of all let us notice that 
\begin{equation*}
\psi^{+}(\V,H) \subseteq H
\end{equation*}
since $\psi^{+}(\V,H)$ is orthogonal to $\V$, as it follows from the fact that $\V$ is null. Let $\V_0:=\psi^{+}(H,H) \subseteq \V$ and let $\V_1$ be the orthogonal complement of $\V_0$ in $\V$.  From the definition of 
$\V_1$ it follows that $\psi^{+}(H,\V_1)$ is orthogonal to $H$ thus $\psi^{+}(H,\V_1)=0$. But $\psi^{+}(\V,\V_1)=0$ as well since $\V$ is null, in other words $\psi^{+}(\V_1, \cdot)=0$ and we conclude that 
$\V_1=0$ using that $\psi^{+}$ is nondegenerate. This proves (i).\\
(ii) is proved by an similar to that in (i), which we leave to the reader. 
\end{proof}
We end this section by the following:
\begin{defi} \label{prod} An SNK-holonomy system $(V^{2m},g,J, \psi^{+},R)$ is said to split as 
$$ V=V_1 \oplus V_2
$$
if $V$ admits an $\mathfrak{h}$-invariant, orthogonal and $J$-invariant splitting $V=V_1 \oplus V_2$ such that 
$$ \psi^{+} \ \mbox{belongs to } \ \lambda^3(V_1) \oplus \lambda^3(V_2). $$
\end{defi}
It is clear that if a SNK-holonomy system splits as $V=V_1 \oplus V_2$ then each of 
$(V_k, g_{\vert V_k}, J_{\vert V_k}, \psi^{+}_{\vert V_k}), k=1,2$ is again a SNK-holonomy system. Note  that factors of dimension $\le 4$ are not permitted by this 
definition since in dimension $2$ or $4$ there are no non-zero holomorphic $3$-forms. Also note that this type of decomposition of a SNK-holonomy system corresponds exactly 
to local products of SNK-manifolds.
%We end this section by the following result which compares isotropic subspaces to null ones in a given SNK-h 
\subsection{The structure of the form $\psi^{+}$}
We are now ready to have a look at the structure of the form $\psi^{+}$ when a complex holonomy reduction is given. The starting point of our approach to the classification 
problem of reducible SNK-holonomy systems is:
% We are interested here 
%in the structure of invariant subspaces of the metric representation $(\g,V)$ and more particularly in those introduced by the following.
\begin{pro} \label{invnk}
Let $\V$ be a proper, $J$-invariant subspace of $(\mathfrak{h},V)$. If $H$ denotes the orthogonal complement of $\V$ in $V$, the following hold:
\begin{itemize}
\item[(i)] $(\psi^{+}_x \circ \psi^{+}_v)w=0$ for all $x$ in $H$ and $v,w$ in $\V$;
\item[(ii)] $(\psi^{+}_x \circ \psi^{+}_y)z$ belongs to $H$ whenever $x,y,z$ are in $H$;
\item[(iii)]  $(\psi^{+}_v \circ \psi^{+}_w)x$ belongs to $H$ for all $x$ in $H$ and $v,w$ in $\V$;
\item[(iv)] $(\psi^{+}_x \circ \psi^{+}_y)v$ is in $\V$, for all $x,y$ in $H$ and all $v$ in $\V$.
\end{itemize} 
\end{pro}
\begin{proof}
Directly from the first Bianchi identity for $R$, combined with the fact that $R(\V,H)=0$ we get 
\begin{equation} \label{curvVH}
R(x,y,v,w)=\langle [\psi^{+}_v, \psi^{+}_w]x,y \rangle-\langle \psi^{+}_vw, \psi^{+}_x,y\rangle
\end{equation}
for all $x$ in $H, y$ in $V$ and $v,w$ in $\V$. All claims are now easy consequences of the fact that $R$ belongs to $S^2(\lambda^{1,1})$ and $\psi^{+}$ is in 
$\lambda^3$, see \cite{Nagy2} for details.
\end{proof}
We can show that any complex reducible SNK-holonomy system contains, up to products,  a null invariant sub-space.
\begin{pro} Let $(V^{2m},g,\psi^{+},R)$ be a complex reducible SNK-holonomy system. Then $V$ splits as 
$$ V=V_1 \oplus V_2
$$
where $V_2$ contains a null invariant space.  
\end{pro}
\begin{proof}
The proof is completed in two steps we shall outline below.\\
{\bf{Step 1}}: Existence of an isotropic invariant subspace.\\
The reducibility of $(\mathfrak{h},V)$ implies the existence of an invariant splitting 
$$ V=E \oplus F
$$
which is moreover orthogonal and stable under $J$. Let $F_0$ be the subspace of $F$ spanned by $\{(\psi^{+}_vw)_F : v,w \ \mbox{in} \ E\}$, where the subscript indicates 
orthogonal projection.  Using Proposition \ref{invnk}, (i) we get 
\begin{equation*}
\psi^{+}_x \psi^{+}_y ((\psi^{+}_vw)_F)+\psi^{+}_x \psi^{+}_y ((\psi^{+}_vw)_E)=0
\end{equation*}
for all $x,y$ in $F$ and whenever $v,w$ are in $E$. But the first summand is in in $E$ by Proposition \ref{invnk}, (iv) while the second 
is in $F$ by (ii) of the same Proposition. therefore both summands vanish individually, and a positivity argument yields then 
$$ \psi^{+}(F,F_0)=0.
$$
In particular $F_0$ is null and from \eqref{mm1} and the $\mathfrak{h}$-invariance of $E$ and $F$ we also get that $F_0$ is $\mathfrak{h}$-invariant. Now $F_0 \neq V$ since 
$\psi^{+} \neq 0$ and if $F_0=0$ we have that $E$ is isotropic, that is $\psi^{+}(E, E) \subseteq E$. \\
{\bf{Step 2 }}: Existence of an invariant nullspace.\\
Using Step 1, we can find a $\mathfrak{h}$-invariant subspace , say $\V$ in $V$, which is isotropic in the sense that 
$\psi^{+}(\V,\V) \subseteq \V$. Let us consider now the $\mathfrak{h}$-invariant tensor $r_1 : \V \to \V$ given by 
$$ \langle r_1v,w\rangle=Tr_{\V} (\psi^{+}_v \circ \psi^{+}_w)
$$
for all $v,w$ in $\V$. Obviously, this is symmetric and $J$-invariant. It follows that there is an orthogonal, $\mathfrak{h}$-invariant splitting 
$$ \V=\V_0 \oplus \V_1
$$
where $\V_0=Ker(r_1)$. Since $\V$ is isotropic, $r_1$ is given by $r_1=-\sum \limits_{k=1}^{d} (\psi_{v_k}^{+})^2v$ for all $v$ in $\V$, where $d:=dim_{\mathbb{R}}\V$. Then Proposition \ref{invnk}, (i) 
implies that $\psi^{+}_x(r_1v)=0$ for all $x$ in $H$ and for all $v$ in $\V$, in other words 
$$ \psi^{+}(H, \V_1)=0.
$$
But the definition of $\V_0$ yields that $\psi^{+}(\V, \V_0)=0$, in particular $\V_0$ is null. We now form the $\mathfrak{h}$ and $J$-invariant subspace $H_1=\V_0 \oplus H_0$ which is easily seen to 
satisfy that 
$$ \psi^{+}(H_1, H_1) \subseteq H_1, \ \psi^{+}(\V_1, \V_1) \subseteq \V_1, \ \psi^{+}(\V_1, H_1)=0. 
$$
In other words $V=\V_1 \oplus H_1$ is a splitting of our holonomy system in the sense of Definition \ref{prod} and the result is proved.
\end{proof}
This can be furthermore (see \cite{Nagy2} for details) refined to:
\begin{pro} \label{spec}
Let $(V^{2m},g,J, \psi^{+},R)$ be an SNK-holonomy system. If it contains a proper invariant nullspace, it splits as 
$$ V=V_1 \oplus V_2
$$
where $V_2$ contains a special invariant subspace.
\end{pro}
Summarising the results obtained up to now we obtain, after an easy induction argument on the irreducible components of $(\mathfrak{h},V)$:
\begin{teo} \label{pdec1} Let $(V^{2m},g,J, \psi^{+},R)$ be a complex reducible SNK-holonomy system. Then $V$ is product of SNK-holonomy systems belonging to one of the following classes:
\begin{itemize}
\item[(i)] irreducible SNK-holonomy systems;
\item[(ii)] SNK-holonomy systems which contain a special invariant subspace.
\end{itemize}
\end{teo}
To advance with the classification of our holonomy systems we need therefore only to discuss the second class present in the Theorem above. We first observe that given a special 
invariant subspace in some SNK-holonomy system one can explicitly determine the curvature along the special subspace.
\begin{pro} \label{curvspec}
Let $(V^{2m},g,J, \psi^{+},R)$ be an SNK-holonomy system containing a special invariant subspace $\V$. Then:
\begin{equation*}
R(\psi^{+}_xy,v_1, v_2, v_3)=\langle [\psi^{+}_{v_1}, [\psi^{+}_{v_2}, \psi^{+}_{v_3}]]x, Jy\rangle
\end{equation*}
whenever $x,y$ are in $H=\V^{\perp}$ and for all $v_k, 1 \le k \le 3$ in $\V$. 
\end{pro}
Actually, Proposition \ref{curvspec} singles out a second Lie algebra of relevance for us, obtained as follows. Consider the subspace $\mathfrak{p} \subseteq \lambda^2(H)$ given by
\begin{equation*}
\mathfrak{p}:=\{\psi^{+}_v : v \ \mbox{in} \ \V \}
\end{equation*}
together with the subspace $\mathfrak{q}$ of $\lambda^{1,1}(H)$ given by $\mathfrak{q}:=[\mathfrak{p}, \mathfrak{p}]$. 
\begin{pro} \label{Lie1}
The following hold:
\begin{itemize}
\item[(i)] $\mathfrak{p} \cap \mathfrak{q}=0$
\item[(ii)] $\mathfrak{r}=\mathfrak{p} \oplus \mathfrak{q}$ is a Lie subalgebra of $\Lambda^2(H)$.
\end{itemize}
\end{pro}
\begin{proof}
From Proposition \ref{curvspec} we get 
$$ [\psi^{+}_{v_1}, [\psi^{+}_{v_2}, \psi^{+}_{v_3}]=\psi^{+}_{R(v_1,v_2)v_3}
$$
for all $v_k, 1 \le k \le 3$ in $\V$. That is $\mathfrak{p}$ is a Lie triple system in the sense that $[\mathfrak{p},[\mathfrak{p},\mathfrak{p}]] \subseteq \mathfrak{p}$ hence 
$\mathfrak{p}+\mathfrak{q}$ is a Lie algebra (see \cite{Helg} for details). To prove (i), let us pick $z$ in $\mathfrak{p} \cap \mathfrak{q} $. Then $z=\psi^{+}_v$ for some 
$v$ in $\V$ and moreover $[z, \mathfrak{p}] \subseteq \mathfrak{p}$. In particular $[\psi^{+}_v, \psi^{+}_{Jv}]=\psi^{+}_w$ for some 
$w$ in $\V$, or equivalently $2(\psi^{+}_v)^2=\psi^{+}_{Jw}$. But the left hand side of this equality is symmetric whilst the right hand is skew-symmetric 
which yields easily that $v=0$.
\end{proof}
The Lie algebra $\mathfrak{r}$ is best though of as the holonomy algebra of a symmetric space of compact type. In the realm of NK-geometry this will turn out to be precisely the case.
\subsection{Reduction to a Riemannian holonomy system}
We shall consider in what follows a SNK-holonomy system $(V^{2m},g,J, \psi^{+},R)$ containing a 
special invariant sub-space $\V$, with orthogonal complement to be denoted by $H$. \\
At this moment need a finer notion of irreducibility for a SNK-holonomy system containing an $\mathfrak{h}$-invariant special subspace. 
Let us now define $R^H : \Lambda^2(H) \to \Lambda^2(H)$ by 
\begin{equation*}
R^H(x,y)=R(x,y)+\psi^{+}_{\psi^{+}_xy}
\end{equation*}
for all $x,y$ in $H$. Note this is well defined because $\psi^{+}(H,H) \subseteq \V$ and $\psi^{+}(\V, H)=H$ and also because $H$ is $\mathfrak{h}$-invariant. Moreover, the Bianchi 
identity for $R$ ensures that $R^H$ is an algebraic curvature tensor on $H$, that is $R^H$ belongs to $\K(\mathfrak{so}(H))$. The second Lie subalgebra of $\Lambda^2(H)$ of interest for us is 
\begin{equation*}
\mathfrak{h}^H=Lie \{ R^H(x,y): x,y \ \mbox{in} \ H\}.
\end{equation*}
\begin{defi} \label{RHS}
$(H,g, \mathfrak{h}^H)$ is called the Riemannian holonomy system associated to $(V^{2m},g,J, \psi^{+}, R)$. It is irreducible if the metric representation $(\mathfrak{h}^H,H)$ is irreducible. 
\end{defi}
The main result in this section is to prove that we can reduce, up to products in the sense of Definition \ref{prod}-hence of Riemannian type-, the study of special SNK-holonomy 
to the case when the associated Riemannian holonomy system is irreducible.
\begin{teo} \label{last}
Let $(V^{2m},g,J, \psi^{+},R)$ be a SNK-holonomy system containing an $\mathfrak{h}$-invariant subspace $\V$. Then $V$ splits as 
$$ V=V_1 \oplus \ldots V_q,
$$
a product of special SNK-holonomy systems, such that each factor has irreducible associated Riemannian holonomy system.
\end{teo}
\begin{proof}
Let $\V$ be a special $\mathfrak{h}$-invariant subspace and let us orthogonally decompose 
$$ H=H_1 \oplus H_2
$$
in invariant subspaces for the action of $\mathfrak{h}^H$. As a straightforward consequence of having $R$ in $S^2(\lambda^{1,1})$ we note that 
\begin{equation} \label{rp}
R^H(\psi^{+}_v)=\frac{1}{2} \psi^{+}_{rv}
\end{equation}
on $H$, for all $v$ in $\V$. Since $R^H(H,H,H_1,H_2)=0$ by assumption and since $r$ is invertible on $\V$ it follows that $\psi^{+}_{\V}H_1$ and $H_2$ are orthogonal, hence 
\begin{equation} \label{eqq1}
\psi^{+}_{\V}H_k \subseteq H_k, k=1,2
\end{equation}
given that $\psi^{+}_{\V}H \subseteq H$. It also folllows that $\psi^{+}(H_1, H_2)=0$. We will need now several preliminary steps.\\
{\bf{Step 1}}: $H_k$ are $J$-invariant for $k=1,2$. \\
From \eqref{eqq1} it follows that $r^{\V}(H_k) \subseteq H_k, k=1,2$
where $r^{\V} : H \to H$ is defined by 
\begin{equation*}
r^{\V}x=-\sum \limits_{k=1}^d(\psi^{+}_{v_k})^2x
\end{equation*}
for all $x$ in $H$, where $\{v_k, 1 \le k \le d=dim_{\mathbb{R}}\V\}$ is some orthonormal basis in $\V$. But $r^{\V}$ has no kernel, as an easy consequence of 
the fact that $\psi^{+}$ is non-degenerate and $\V$ is special. It follows that $r^{\V}_{\vert H_k} : H_k \to H_k$ is injective, hence surjective, that is 
$$ r^{\V}(H_k)=H_k
$$
for $k=1,2$. Again from \eqref{eqq1} we find that $r^{\V}(JH_k) \subseteq H_k$ hence after applying $(r^{\V})^{-1}$ (which as we have seen preserves $H_k$) we get that 
$JH_k \subseteq H_k$ for $k=1,2$. \\
{\bf{Step 2}}: A first decomposition of $V$\\
Let $\V^{\prime}:=\V \oplus H_1$ and consider the orthogonal and $J$-invariant splitting 
$$ V=\V^{\prime} \oplus H_2.
$$
Since $R^H(H_1,H_2,H,H)=0$ and also 
because $\psi^{+}(H_1,H_2)=0$ we find that \\
$R(H_1,H_2,H,H)=0$. In other words we have $R(H,H)H_2 \subseteq H_2$ and also $R(H,\V)H_2=0$ since $\V$ is $\mathfrak{h}$-invariant and $R$ is symmetric.
Now the Bianchi identity for $R$, combined again with the $\mathfrak{h}$-invariance of $\V$ and its nullity gives 
$$ R(v,w)x=[\psi^{+}_v, \psi^{+}_w]x
$$
for all $v,w$ in $\V$ and for all $x$ in $H$. Using again \eqref{eqq1} we find that $R(\V,\V)H_2 \subseteq H_2$ so $H_2$ is actually $\mathfrak{h}$-invariant. \\
{\bf{Step 3}}: The decomposition of $\V$\\
Thus we may apply Proposition \ref{invnk} (i) to the splitting $V=\V^{\prime} \oplus H_2$ whence 
$$ \psi^{+}_x \psi^{+}_vw=0
$$
for all $x$ in $H_2$ and whenever $v,w$ belong to $\V^{\prime}$.  Let us define $\V_k=\psi^{+}(H_k, H_k) \subseteq \V$ for $k=1,2$ and notice these are $J$-invariant. Then 
$$ \psi^{+}(H_2, \V^1)=0
$$
in particular $\V^1$ and $\V^2$ are orthogonal, after taking the scalar product with elements in $H_2$.  
and since $H_1$ and $H_2$ play equal roles we also get that $\psi^{+}$. Therefore 
$$ \V=\V^1 \oplus \V^2
$$
after using that $\psi^{+}(H, H)=\V$ and $\psi^{+}(H_1,H_2)=0$.\\
Now since $H_1$ and $H_2$ play equal roles, by repeating the arguments above we also get that $\psi^{+}(\V^2,H_1)=0$. \\
{\bf{Step 4}}: Proof of the Theorem.\\
We consider the orthogonal and $J$-invariant splitting 
$$V=V_1 \oplus V_2 
$$ 
where $V_k:=\V^k \oplus H_k, k=1,2$. From Step 3 it follows that $\psi^{+}$ belongs to 
$$\lambda^3(V_1) \oplus \lambda^3(V_2)$$ 
and moreover that $H_k$ are invariant under  $\mathfrak{h}$ for $k=1,2$. Let us prove that, say $\V^1$, is $\mathfrak{h}$-invariant too. From \eqref{mm1} we get  
$$ R(x,y)(\psi^{+}_{x_1}x_2)=\psi^{+}_{R(x,y)x_1}x_2+\psi^{+}_{x_1}R(x,y)x_2
$$
for all $x,y$ in $V$ and $x_1,x_2$ in $H_1$ and for all $x,y$ in $H$. Since $H$ is $\mathfrak{h}$-invariant and we have seen that $\psi^{+}(H,H_1)=\psi^{+}(H_1,H_1)=\V^1$, the invariance of 
of $\V^1$ follows and that of $\V^2$ is proved analogously. We conclude that $V^k$ are invariant under $\mathfrak{h}$, for $k=1,2$.\\
The claim follows now by induction on the number of irreducible components of $(\mathfrak{h}^H,H)$, provided we show that $\mathfrak{h}^H$ does not fix any vector in $H$. But this is implied by 
\eqref{rp}, for such a vector, say $x_0$ in $H$ would satisfy then $\psi^{+}(\V,x_0)=0$. Taking scalar products with elements in $H$ this yields $\psi^{+}(x_0,H)=0$ as $\psi^{+}(H,H) \subseteq \V$.  This is to say that 
$\psi^{+}_{x_0}=0$ hence $x_0=0$ since $\psi^{+}$ is non-degenerate. This proves the absence of fixed points for the representation $(\mathfrak{h}^H,H)$ and the proof is finished.
\end{proof}
More properties of the non-Riemannian holonomy representation $(\mathfrak{h},V)$ can be now formulated.
\begin{pro} \label{iredfiber}
Let $(V^{2m},g,J,\psi^{+})$ be a special SNK-holonomy system containing an invariant subspace $\V$ and such that the Riemannian holonomy system $(\mathfrak{h}^H,H)$ is irreducible. Then 
the representation $(\mathfrak{h},\V)$ is irreducible over $\mathbb{C}$. 
\end{pro}
\begin{proof}
Let us suppose that $\V$ is not irreducible and split $\V=\V^1 \oplus \V^2$ as orthogonal sum of $J$-stable and $\mathfrak{h}$-invariant subspaces. Then $R(\psi^{+}_xy,v,v_1, v_2)=0$ for all 
$x,y$ in $H, v$ in $\V$ and for all $v_k$ in $\V^k, k=1,2$. Using Proposition \ref{curvspec} this leads to 
$$ [\psi^{+}_v, [\psi^{+}_{v_1}, \psi^{+}_{v_2}]]=0
$$
whenever $v$ belongs to $\V$ and $v_1, v_2$ are in $\V^1, \V^2$ respectively. An easy invariance argument which can be found in \cite{Nagy2}, page 492 yields 
\begin{equation} \label{temp1}
\psi^{+}_{v_1} \psi^{+}_{v_2}=0
\end{equation}
for all $v_k$ in $\V^k, k=1,2$. We form $H^k:=\psi^{\V^k}H, k=1,2$ which are therefore orthogonal, $J$-invariant and such that $H=H^1 \oplus H^2$. The spaces $H^k, k=1,2$ are actually invariant under $\mathfrak{h}$ 
fact which follows from \eqref{mm1} and the $\mathfrak{h}$-invariance of $\V^k, k=1,2$. Now by using \eqref{temp1} we obtain that 
\begin{equation*}
\psi^{+}_{\psi^{+}_xy} H^k \subseteq H^k
\end{equation*}
for all $x,y$ in $H$, since $ \psi^{+}_{\V^k}H \subseteq H$ and $\psi^{+}(H,H) \subseteq \V$. This means, see also the definition of the tensor $R^H$, that $H^k, k=1,2$ are invariant under 
$\mathfrak{h}^H$ hence they cannot be both proper since $(\mathfrak{h}^H,H)$. If $H^1=0$ for instance, a routine argument leads to $\V^1=0$, a contradiction and the proof is finished.
\end{proof}
\subsection{Metric properties}
Let  $(V,g,J, \psi^{+},R)$ be a SNK-holonomy system with special invariant subspace $\V$ and such that the associated representation $(\mathfrak{h}^H,H)$ is irreducible. We have 
seen that $(\mathfrak{h},\V)$ is an irreducible representation, but let us notice that $(\mathfrak{h},H)$ can be still reducible. Instead of attacking the problem of irreducibility of $(\mathfrak{h},H)$ directly 
we shall use the Ricci tensor of our holonomy system as an indicator in this direction. In the context of nearly-K\"ahler manifolds this results in metric and topological information.\\
Note that the tensors $Ric, \overline{Ric}$ and $C$ are defined here exactly as in the beginning of this section, and are still subject to \eqref{ricf} and \eqref{chern} as it follows from the corresponding proofs. Their invariance under 
$\mathfrak{h}$ is an easy consequence of \eqref{ricf}. We shall prove here that one reduces the discussion to the case when $C$ has only a few eigenvalues but skip most of the technical details.
\begin{pro} \label{eigC} The tensor $C$ has at most $3$-eigenvalues. 
\end{pro}
\begin{proof}
Directly from its definition $C$ preserves the splitting $V=\V \oplus H$. Since $C$ is $\mathfrak{h}$-invariant and $(\mathfrak{h},\V)$ is irreducible it follows that $C_{\vert \V}=\lambda 1_{\V}$ for some 
real number $\lambda$. Using \eqref{chern} we find that the symmetric and $J$-invariant tensor $S:H \to H, S:=C_{\vert H}+\frac{\lambda}{2}$ satisfies 
$$ S(\psi^{+}_vx)=-\psi^{+}_vSx
$$
hence 
\begin{equation}  \label{sqC}
S^2(\psi^{+}_vx)=\psi^{+}_vS^2x 
\end{equation}
for all $x$ in $H$ and for all $v$ in $\V$.  Since $C$ is $\mathfrak{h}$-invariant, so is $S^2$, and in fact the latter turns out to be $\mathfrak{h}^H$-invariant, after using \eqref{sqC}. It follows that $S^2$ is a multiple 
of the identity and the claim follows.
\end{proof}
Let us consider the tensor $r^{\V} : H \to H $ given by 
$$ r^{\V}=-\sum \limits_{k=1}^{2d}(\psi^{+}_{v_k})^2
$$
where $2d=dim_{\mathbb{R}}\V$ and $\{v_k, 1 \le k \le 2d\}$ is some orthonormal basis in $\V$. It is easy to see $r^{\V}$ is positive, $J$-commuting, symmetric and $\mathfrak{h}$-invariant 
and moreover that $r_{\vert H}=2r^{\V}$. \par
When $C$ has only one eigenvalue, then $C=0$ and one can show \cite{Nagy2} by using \eqref{ricf} that both $Ric$ and $r$ must actually by diagonal on $V$. \par
It remains to investigate the cases when $C$ has two, respectively three eigenvalues which are described below.
\begin{pro} \label{eCfin}
The following hold:
\begin{itemize}
\item[(i)] If $C$ has exactly two eigenvalues then there exists $k >0$ such that  $r^{\V}=k 1_H$ and moreover the eigenvalues together with the corresponding eigenspaces for the tensors 
$r,C, Ric$ are given in the following table:
%\begin{equation*}
$$ \begin{tabular}{ | r |c | c | c |c |} 
\hline
Eigenvalue  &  r &  Ric &  C & Eigenspace \\
\hline 
$\lambda_1 $ & $\frac{n-d}{d}k$&  $\frac{n+7d}{4d} k $&  $\frac{4(n-3d)}{d}k$ & $\V$ \\
\hline
$\lambda_2 $ & $2k$ & $\frac{n+2d}{2d}k$& $-\frac{2(n-3d)}{d}k$ & $H$ \\
\hline
\end{tabular}$$
%\end{equation*}
\item[(ii)] If $C$ has exactly $3$-eigenvalues we have a $\mathfrak{h}$-invariant, orthogonal and $J$-invariant splitting 
\begin{equation*}
V=\V \oplus H_1 \oplus H_2
\end{equation*} 
such that $\psi^{+}(H_1,H_2)=\psi^{+}(H_2,H_2)=0$ and $\psi^{+}(H_1,H_2)=\V$.
\item[(iii)] in all cases, including when $C=0$,  we have 
\begin{equation*}
Ric^H=\mu g_{\vert H}
\end{equation*}
for some $\mu >0$, where $Ric^H$ denotes the Ricci contraction of $R^H$. 
\end{itemize}
\end{pro}
We refer the reader to \cite{Nagy2} for details of the proof and also mention that in case (ii) above it is also possible to display the relations between the eigenvalues of all relevant symmetric endomorphisms. 
\subsection{The final classification}
By Theorem \ref{last} it is enough to consider a SNK-holonomy system $(V^{2m},g,J, \psi^{+},R)$ with a special invariant subspace $\V$ and such that the Riemannian holonomy system 
$(\mathfrak{h}^H,H)$ is irreducible, in the notations of the previous section.  
We shall combine here some well known results for irreducible Riemannian holonomy systems and the information which is available in our present context. To progress in this direction we will
first make clear the relationship the Lie algebra $\mathfrak{r}$ and the Riemannian holonomy algebra $\mathfrak{h}^H$. 
\begin{pro} \label{ideal+}The following hold:
\begin{itemize}
\item[(i)] $\mathfrak{r}$ is an ideal in $\mathfrak{h}^H$;
\item[(ii)] $R^H=\frac{\nu}{2} 1_{\mathfrak{r}}$ on $\mathfrak{r}$ for some $\nu>0$ such that $r_{\vert \V}=\nu 1_{\V}$.
\end{itemize}
\end{pro}
\begin{proof}
(i) That $\mathfrak{p} \subseteq \mathfrak{h}^H$ follows from \eqref{rp} hence $\mathfrak{q}=[\mathfrak{p}, \mathfrak{p}]$ is contained in $\mathfrak{h}^H$ as well. To see that 
$\mathfrak{r}$ is an ideal we use \eqref{mm1} to observe that 
$$ [R(x,y), \psi^{+}_v]
$$
belongs to $\mathfrak{p}$ for all $x,y$ in $H$ and for all $v$ in $\V$. The definition of $R^H$ and that of $\mathfrak{r}$ now yields
$$ [R^H(x,y), \psi^{+}_v]
$$
in $\mathfrak{r}$ and it easy to conclude by using the Jacobi identity and Proposition \ref{Lie1}.\\
(ii) First of all let us observe that $r$ preserves $\V$, since the latter is special. The $\mathfrak{h}$-invariance of $r$, together with the irreducibility of  $(\mathfrak{h},\V)$ implies the existence 
of a constant $\nu>0$ such that $r_{\vert \V}=\nu 1_{\V}$. For the action of $R^H$ on $\mathfrak{p}$ our claim follows now from \eqref{rp}. 
%Let us prove a slightly more general fact. 
To prove it on $\mathfrak{q}$ we observe we proceed as follows. Given that $R$ belongs to $S^2(\Lambda^2)$ we can alternatively rewrite \eqref{mm1}
under the form 
\begin{equation*}
R(\psi^{+}_{x_1}x_2,x_3)+R(x_2, \psi^{+}_{x_1}x_3)=R(x_1, \psi^{+}_{x_2}x_3)
\end{equation*}
whenever $x_k$ are in $V, k=1,2,3$. Using this for $x_1=v, x_2=\psi^{+}_we_i, x_3=e_i$ where $v,w$ are in $\V$ and $\{e_i \}$ is some orthonormal basis in $H$ yields after summation and rearrangement of terms  
$$ \sum \limits_{i} R(e_i, [\psi^{+}_v, \psi^{+}_w] e_i)=R(v,rw)
$$
that is $R([\psi^{+}_v, \psi^{+}_w])=\frac{\nu}{2}[\psi^{+}_v, \psi^{+}_w]$ on $H$, after also using \eqref{curvVH}. After a straightforward invariance argument this 
leads to $R^H([\psi^{+}_v, \psi^{+}_w])=\frac{\nu}{2}[\psi^{+}_v, \psi^{+}_w]$ and the claim is proved.
\end{proof}
Now we recall that for a given metric representation $(\mathfrak{l},W)$ of some Lie algebra $\mathfrak{l}$ we have the Berger list of irreducible Riemannian holonomies, 
$$ \begin{tabular}{ | r |c | c | } 
\hline
$\mathfrak{l}$  &  W &  \\ 
\hline 
$\mathfrak{so}(n) $ & $\mathbb{R}^n$ & \\
$\mathfrak{u}(m) $ & $\mathbb{R}^{2m}$ & \\
$\mathfrak{su}(m)$ & $\mathbb{R}^{2m}$ &\\
$\mathfrak{sp}(m) \oplus \mathfrak{sp}(1)$ & $\mathbb{R}^{4m}$ & \\
$\mathfrak{sp}(m)$ & $\mathbb{R}^{4m}$ & \\
$\mathfrak{spin}(7)$ & $\mathbb{R}^8$ & \\
$ \mathfrak{g}_2 $ & $\mathbb{R}^7$ & \\
\hline
\end{tabular}$$
and proceed to the final classification result, hence proving Theorem \ref{intromain} in the introduction.
\begin{teo} \label{finclassnk}
Let $(M^{2m},g,J)$ be a complete SNK-manifold. Then $M$ is, up to finite cover, a Riemannian product whose factors belong to the following classes:
\begin{itemize}
\item[(i)] homogeneous SNK-manifolds;
\item[(ii)] twistor spaces over positive quaternionic K\"ahler manifolds;
\item[(iii)] $6$-dimensional $SNK$-manifolds.
\end{itemize}
\end{teo}
\begin{proof}
First of all, because  $M$ must be compact with finite fundamental group we may assume up to a finite cover that it is simply connected. Pick now a point $x$ in $M$ and consider the holonomy representation 
$(Hol_x(\bnabla),T_xM)$. If this is irreducible over $\mathbb{R}$ or complex irreducible but real reducible we conclude by  Theorem \ref{irednk}.\\
If the holonomy representation of $\bnabla$ is reducible over $\mathbb{C}$ at $x$ we consider the SNK-holonomy system $(T_xM, g_x, J_x, \psi^{+}_x, \overline{R}_x)$. 
By the theorem of Ambrose-Singer we have that $\mathfrak{h} \subseteq \mathfrak{hol}_x(\bnabla)$. Also note that when {\it{starting}} the whole reduction procedure leading to Theorem \ref{last} from 
some $\mathfrak{hol}_x(\bnabla)$ complex invariant sub-space we end up with a splitting 
$$ T_xM=V_1 \oplus \ldots V_p
$$ 
with the properties in Theorem \ref{last} and which is furthermore $\mathfrak{hol}_x(\bnabla)$-invariant together with the special sub-spaces contained in each factor. Using parallel transport w.r.t to $\bnabla$ this extends 
to a $\bnabla$-parallel decomposition, which is actually $\nabla$-parallel since the torsion of $\bnabla$ is split along along the decomposition. The splitting theorem of deRham applies and 
our study is reduced to that of SNK-manifolds such that at each point, we have an invariant special subspace having the property that the associated holonomy system is irreducible in the sense of Definition
\ref{RHS}. Then we show (see \cite{Nagy2}) that $M$ is the total space of a Riemannian submersion 
\begin{equation*}
\pi : (M,g) \to (N^{2n},h)
\end{equation*}
whose fibres are totally geodesic and w.r.t the induced metric and almost complex structure are simply connected, irreducible and compact Hermitian symmetric spaces. The tensor $R^H$ defined in section 5.3 projects 
onto the Riemann curvature tensor of $h$ making of $\mathfrak{h}^H$ a subalgebra of $\mathfrak{hol}(N,h)$ at each point. Using this one shows that the Riemannian manifold $(N^{2n},h)$ is irreducible 
and Proposition \ref{eCfin}, (iii) gives that $h$ is Einstein of positive scalar curvature.  \\
If the tensor $C$ has three eigenvalues a short manipulation of the second Bianchi identity for $\bnabla$ combined with the algebraic structure of the splitting of $TM$ given in (ii) of Proposition \ref{eCfin} gives that 
$\bnabla$ is an Ambrose-Singer connection.  \\
Let us discuss now the case when $C$ has two eigenvalues or it vanishes. \\
If $(N^{2n},h)$ is a symmetric space the relations of O'Neill (see \cite{Besse}) for the Riemannian submersion $\pi : M \to N$ essentially say that the curvature tensor $\overline{R}$ is an explicit 
expression of $\psi^{+}, R^H$ and the curvature of the fibre (which is a symmetric space) along the $\bnabla$-parallel splitting $TM=\V \oplus H$. This information results easily in having $\bnabla$ an Ambrose-Singer connection 
again. \\
Suppose now that $(N^{2n},h)$ is not a symmetric space. Then the holonomy representation of $(N,h)$, at the Lie algebra level, corresponds to the representation $(\mathfrak{h},H)$ and as it is well known (see \cite{redbook}), 
must be one of the entries in the Berger list above. We shall mainly use now Proposition \ref{ideal+}. The Lie algebras $\mathfrak{su}(m),\mathfrak{sp}(m),\mathfrak{spin}(7),\mathfrak{g}_2$ are excluded 
because their curvature tensors are Ricci-flat. 
Supposing that $\mathfrak{h}^H=\mathfrak{u}(n)$, given that the ideal $\mathfrak{r}$ is at least two-dimensional we can only have $\mathfrak{r}=\mathfrak{su}(n), \mathfrak{u}(n)$. But using Proposition \eqref{ideal+} combined with 
the fact that a curvature tensor of K\"ahler-Einstein type is completely determined by its restriction to $\mathfrak{su}(m)$ we find that $(N^{2n},h)$ is symmetric a contradiction. Similar arguments together with results 
from \cite{BeMor} enable us to conclude in the case when $\mathfrak{h}^H=\mathfrak{so}(2n)$ or when $\mathfrak{r}=\mathfrak{sp}(m)$. The only case remaining is when $\mathfrak{h}^H=\mathfrak{sp}(m) \oplus 
\mathfrak{sp}(1)$ and $\mathfrak{r}=\mathfrak{sp}(1)$ which implies that  $\V$ is of rank $2$. Hence the holonomy group of $\bnabla$ is contained in $U(1) \times U(2m)$ and the fact that $(M^{2m},g,J)$ is a twistor space 
over a positive quaternionic-K\"ahler manifold has been proved in \cite{Nagy1}.
\end{proof}
\section{Concluding remarks}
Due to considerations of space and time we have omitted to present some important facts related to nearly-K\"ahler geometry. The first aspect is related to the 
construction of examples, which is implicit in our present treatement. Given a quaternion-K\"ahler manifold $(M^{4m},g,Q)$ of positive scalar curvature the Salamon twistor construction \cite{SMS1} yields a K\"ahler manifold $(Z,h,J)$ such that there is a Riemannian submersion with totally geodesic, complex fibres 
$$ S^2 \hookrightarrow (Z,h) \to (M,g).$$
Using this structure and a canonical variation of the metric $h$ it has been shown \cite{SElls} that $Z$ admits an SNK-metric. Enlarging the context to that of 
Riemannian submersions with complex, totally geodesic fibres from K\"ahler manifolds one draws a similar conclusion \cite{Nagy1} (see also \cite{Nagy3} for some 
related facts). This observation can also be used to describe the homogeneous SNK-manifolds which appear in Theorem \ref{finclassnk} as twistor spaces 
in the sense of \cite{Rawnsley} over symmetric spaces of compact type. 

We have also omitted to discuss NK-structures in dimension $6$, which have a rich geometry, although not yet fully understood. If $(M^6,g,J)$ is a strict 
NK manifold the structure group automatically reduces to $SU(3)$, hence $c_1(M,J)=0$ \cite{g3}. Moreover the metric $g$ must be Einstein \cite{g3}, of positive 
scalar curvature. There is also a $1:1$ correspondence between $6$-dimensional manifolds admitting rela Killing spinors and SNK-structures in dimension $6$ 
\cite{grunewald}. 

In the compact case the only 
known examples are homogeneous, and various characterisations of these instances are available \cite{But1,MNS1}. 
However, compact examples with two conical singularities have been very recently constructed in \cite{FernMunIva} building on the fact that one can suitably 
rotate the $SU(2)$-structure of a Sasakian-Einstein $5$-manifold and then extend it to a $SU(3)$-structure of NK-type through a generalised cone construction. 
$\\$
$\\$
{\bf{Acknowledgments}}: 
%I thank Andrei Moroianu and Uwe Semmelmann for colaboration on some of the results presented here. 
Part of this material was elaborated during the "77\`eme rencontre entre physiciens th\'eoriciens et math\'ematiciens" held in Strasbourg in 2005. Warm thanks to the organisers and especially to V.Cort\`es for a beautiful 
and inspiring conference. I am also grateful to V.Cort\`es and A.Moroianu for useful comments.

\end{document}